\documentclass[11pt,leqno,draft]
{article}
\input epsf.tex
\usepackage{amsthm,array,
amsmath, amssymb, amscd, amsfonts,
latexsym, verbatim}

\input xy
\xyoption{all}

\usepackage{euscript}

\begin{document}
\input{amssym.def}

%%%%%%%%%%%%%%%%%%%%%%%%%%%%%%%%%%%%%%%%%

\newtheorem{theorem}{Theorem}
\newtheorem{corollary}{'orollary}
\newtheorem{proposition}{Proposition}
\newtheorem{lemma}{Lemma}
\newtheorem{maintheorem}{'heorem}
\newtheorem{maincorollary}{Corollary}
\newtheorem{remark}{Remark}
\newtheorem{definition}{Definition}
\newtheorem{example}{\bf Example}

%%%%%%%%%%%%%%%%%%%%%%%%%%%%%%%%%

\newcommand{\C}{\mathbb C}
\newcommand{\Ha}{\mathbb H}
\newcommand{\R}{\mathbb R}
\newcommand{\OO}{\mathbb O}
\newcommand{\g}{\mathfrak g}
\newcommand{\h}{\mathfrak h}
\newcommand{\z}{\mathfrak z}
\newcommand{\ssl}{\mathfrak {sl}}
\newcommand{\spl}{\mathfrak {sp}}
\newcommand{\su}{\mathfrak {su}}
\newcommand{\so}{\mathfrak {so}}
\newcommand{\gl}{\mathfrak {gl}}
\newcommand{\GL}{\rm {GL}}
\newcommand{\Or}{\rm {O}}
\newcommand{\U}{\rm {U}}
\newcommand{\Sp}{\rm {Sp}}
\newcommand{\Id}{\rm {Id}}
\newcommand{\I}{\rm {I}}
\newcommand{\ka}{\mathfrak k}
\newcommand{\te}{\mathfrak t}
\newcommand{\un}{\mathfrak u}
\newcommand{\Mat}{\rm {Mat}}
\newcommand{\ov}{\overline}
\newcommand{\p}{\mathfrak p}
\newcommand{\ru}{{\rm rad}_u}
\newcommand{\rd}{\rm rad}
\newcommand{\al}{\mathfrak a}
\newcommand{\pt}{\!\cdot\!}
\newcommand{\pe}{\mathfrak p}
\newcommand{\gR}{{\g}^{}_{\R}}
\newcommand{\kR}{{\ka}^{}_{\R}}
\newcommand{\pR}{{\pe}^{}_{\R}}
\newcommand{\GRc}{{\rm Ad}(\g^{}_{\R})}
\newcommand{\GR}{{G}\hskip -.1mm_{\R}}
\newcommand{\cNv}{{\cal N}(V)}
\newcommand{\cNg}{{\cal N}(\g)}
\newcommand{\cNp}{{\cal N}(\p)}
\newcommand{\cNgr}{{\cal N}(\gR)}

%%%%%%%%%%%%%%%%%%%%%%%%%%%%%%%
\newcommand{\blist}{
\begin{list}{}
{ \setlength{\topsep}{2.5mm}
\setlength{\parsep}{-1mm}
\setlength{\itemindent}{0mm}
\setlength{\labelwidth}{10mm}
\setlength{\leftmargin}{15mm}
\setlength{\labelsep}{2mm} } }
%%%%%%%%%%%%%%%%%%%%%%%%%%%%%%%%%%

\title{\vskip -10mm ${\ }$ \vskip -2cm\bfseries Self-dual  projective
algebraic varieties\\
associated with symmetric spaces{\tiny
$^{^{\clubsuit}}$}\footnotetext{\hskip -2mm{\tiny
$^{^{\clubsuit}}$}Published in: {\it Algebraic
Transformation Groups and Algebraic Varieties}, Enc.
Math. Sci., Vol. 132, Sub\-se\-ries {\it Invariant
Theory and Algebraic Transformation Groups}, Vol. III,
Springer Verlag, 2004, 131--167.}}

\author{
{\bf Vladimir~L.~Popov}\footnote {Partly supported by
grant \hskip 2mm${}_{-}$\hskip -1mm${}_{-}$\hskip
-5.95mm
%%{\rus
HIII-123.2003.01 (Russia), ESI (Vienna,
Austria), and
ETH (Z\"urich, Swit\-zerland).}\\[2.5pt]
{\small Steklov Mathematical
Institute,}\\[-2.5pt]
{\small Russian Academy of
Sciences,}\\[-2.5pt]
{\small Gubkina 8,  Moscow 117966,
Russia }\\[-2.2pt]
{\small {\tt popovvl@orc.ru}}
\\
\\
{\bf Evgueni~A.~Tevelev}\footnote {Partly
supported
by ESI (Vienna, Austria).}\\[2.5pt]
{\small Department of Mathematics,}\\[-2.5pt]
{\small The University of Texas,}\\[-2.5pt]
{\small Austin, Texas 78712, USA }\\[-2.2pt]
{\small {\tt tevelev@mail.ma.utexas.edu}}
}
\date{October 31, 2003}
\maketitle

%%\tableofcontents

%%\pagestyle{myheadings} \markboth {
%%V.~L.~POPOV, E.~A.~TEVELEV} {SELF-DUAL
%%ALGEBRAIC VARIETIES}

\pagestyle{myheadings} \markboth{\hskip
-.5cm\centerline{\rm V.~L.~POPOV,
E.~A.~TEVELEV} \hskip -.2cm}{\hskip
-.1cm\centerline{\rm SELF-DUAL ALGEBRAIC
VARIETIES}\hskip -.5cm}

%%\

%%\vskip -2cm

%%\

\begin{abstract}
%\vskip -115mm
%{\hfill To appear in: Vol. 3 of Springer
%Verlag
%%{\it Interesting
%%Algebraic Varieties}}
%%
%%{\hfill \it Arising in Algebraic
%%Transformation Groups Theory},
%%{\hfill
%Enc. Math. Sci. Subser.}
%
%{\hfill {\it Invariant Theory and
%Algebraic Transformation Groups}, 2004} %%
%%{\hfill Vol. 3, 2004}
%\vskip 110mm
We
discover a class of projective self-dual
algebraic varieties. Namely, we consider
actions of isotropy groups of complex
symmetric spaces on the projectivized
nilpotent varieties of isotropy modules.
For them, we classify all orbit closures
$X$ such that $X=\check X$ where $\check
X$ is the projective dual of $X$. We give
algebraic criteria of projective
self-duality for the considered orbit
closures.
\end{abstract}

\section*{1. Introduction}
\label{introduction}  Under different
guises dual varieties of projective
algebraic varieties have been
consi\-de\-red in various branches of
mathematics for over a hundred years. In
fact, the dual variety is the
generalization to algebraic geometry of
the Legendre transform in classical
mechanics, and the Biduality Theorem
essentially rephrases the duality between
the Lagrange and Hamilton--Jacobi
approaches in classical mechanics.

Let $X$ be an $n$-dimensional projective
subvariety of an $N$-dimensi\-onal
projective space $\mathbf P$, and let
$\check{X}$ be the dual variety of $X$ in
the dual projective space $\check{\mathbf
P}$. Since various kinds of geometrically
meaningful unusual behavior of hyperplane
sections are manifested more explicitly in
terms of dual varieties, it makes sense to
consider their natural invariants. The
simplest invariant of $\check X$ is its
dimension $\check n$. ``Typically'',
$\check n=N-1$, i.e., $\check X$ is a
hypersurface. The deviation from the
``typical" behavior admits a geometric
interpretation: if $\check X$ is not a
hypersurface and, say,
$\operatorname{codim}\check X = s+1$, then
$X$ is uniruled by $s$-planes.

Assume that $X$ is a smooth variety not
contained in a hyperplane. Then
$n\leqslant\check n$, by \cite[Chapter
1]{Z}.  If the extremal case $n=\check n$
holds and $n=\check n \leqslant 2N/3$,
then, by \cite{Ei1}, \cite{Ei2}, such $X$
are classified by the following list:

\blist \item[(i)] hypersurfaces,
\item[(ii)] $\mathbf P^1 \times \mathbf
P^{n-1}$ embedded in $\mathbf P^{2n-1}$ by
the Segre embedding, \item[(iii)] the
Grassmannian of lines in $\mathbf P^4$
embedded in $\mathbf P^9$ by the Pl\"ucker
map,~or \item[(iv)] the $10$-dimensional
spinor variety of $4$-dimensional linear
subspaces on a non\-sin\-gular
$8$-dimensional quadric in $\mathbf
P^{15}$.
\end{list}

According to Hartshorne's famous
conjecture, if $n>2N/3$, then $X$ is a
complete intersection, and hence, by
\cite{Ei2}, $\check X$ is a hypersurface.
Therefore it is plausible that the
above-stated list  contains every smooth
$X$ such that $n=\check n$. Furthermore,
in cases (ii), (iii) and (iv), the variety
$X$ is {\it self-dual} in the sense that,
as an embedded variety, $\check X$ is
isomorphic to $X$. In case (i), it is
self-dual if and only if it is a quadric.
Thus, modulo Hartshorne's conjecture, the
above-stated list gives the complete
classification of all {\it smooth}
self-dual projective algebraic varieties.
In particular it shows that there are not
many of them.

In \cite{P2}, it was found a method for
constructing many {\it nonsmooth}
self-dual projective algebraic varieties.
It is related to algebraic transformation
groups theory: the self-dual projective
algebraic varieties appearing in this
construction are certain projectivized
orbit closures of some linear actions of
reductive algebraic groups. We recall this
construction.

Let $G$ be a connected reductive algebraic
group, let $V$ be a finite dimensional
algebraic $G$-module and let $B$ be a
nondegenerate symmetric $G$-invariant
bilinear form on~$V$. We assume that the
ground field is $\mathbb C$. For a subset
$S$ of $V$, put $S^{\perp}:=\{v\in V\mid
B(v, s)=0 \ \forall \ s\in S\}$. We
identify $V$ and $V^*$ by means of $B$ and
denote by $\mathbf P$ the associated
projective space of $V=V^*$. Thereby the
projective dual $\check X$ of a Zariski
closed irreducible subset $X$ of $\mathbf
P$ is a Zariski closed subset of $\mathbf
P$ as well. Let $\cNv$ be the {\it
null-cone} of $V$, i.e.,
$$
\cNv:=\{v\in V\mid 0\in \overline{G\cdot
v}\},
$$
where bar stands for Zariski closure in
$V$.
\begin{theorem} {\rm (\cite[Theorem
1]{P2})}\label{construction} Assume that
there are only finitely many $G$-orbits in
$\cNv$. Let $v\in \cNv$ be a nonzero
vector and let $X:={\mathbf
P}(\overline{G\cdot v})\subseteq \mathbf
P$ be the projectivization of its orbit
closure. Then the following properties are
equivalent: \blist \item[{\rm (i)}] $X=
\check X$. \item[{\rm (ii)}] $({\rm
Lie}(G)\cdot v)^{\perp}\subseteq \cNv$.
\end{list}
\end{theorem}

Among the modules covered by this method,
there are two natu\-ral\-ly allo\-ca\-ted
classes, namely, that of the adjoint
modules and that of the isotropy modules
of symmetric spaces. In fact the first
class is a subclass of the second one;
this subclass has especially nice
geometric properties and is studied in
more details. Projective self-dual
algebraic varieti\-es associated with the
adjoint modules by means of
Theorem~\ref{construction} were explicitly
classified and studied in \cite{P2}. In
particular, according to \cite{P2}, these
varieties are precisely the projectivized
orbit closures of nilpotent elements in
the Lie algebra of $G$
%%${\rm Lie}(G)$
that are
 {\it distinguished} in the sense of Bala
 and Carter (see below Theorem~\ref{dist}).
This purely algebraic notion plays an
important role in the Bala--Carter
classification of nilpotent
elements,~\cite{BaCa}.

 For the isotropy modules of symmetric
 spaces, in \cite{P2} it was introduced
the notion of {\it
$(-1)$-dis\-tin\-gui\-shed element} of a
$\mathbb Z_2$-graded semisimple Lie
algebra and shown that the projective
self-dual  algebraic varieties associated
with such modules by means of
Theorem~\ref{construction} are precisely
the projectivized orbit closures of
$(-1)$-distinguished elements. Thereby
classification of such varieties was
reduced to the problem of classifying
$(-1)$-distinguished elements. In
\cite{P2}, it was announced that the
latter problem will be addressed in a
separate publication.

The goal of the present paper is to give
the announced classification: here we
expli\-cit\-ly classify
$(-1)$-distinguished elements of all
$\mathbb Z_2$-graded complex semisimple
Lie algebras; this yields the
classification of all projective self-dual
algebraic varieti\-es associated with the
isotropy modules of symmetric spaces by
means of the above-stated construction
from \cite{P2}. In the last section we
briefly discuss some geometric properties
of these varieties.

Notice that there are examples of singular
projective self-dual algebraic varieties
constructed in a different way. For
instance, the Kummer surface in $\mathbf
P^3$, see \cite{GH}, or the Coble quartic
hypersurface in $\mathbf P^7$, see
\cite{Pa}, are projective self-dual. A
series of examples is given by means of
``Pyasetskii pairing'', \cite{T}
(e.g.,\;projectivization of the cone of
$n\times m$-matrices of rank $\leqslant
n/2$ is projective self-dual for
$n\leqslant m$, $n$ even).

Given a projective variety, in general it
may be difficult to explicitly identify
its dual variety. So our classification
contributes to the problem of finding
projective varieties $X$ for which $\check
X$ can be explicitly identified.
%% providing
%%%many examples of such $X$.
Another application concerns the problem
of explicit describing the projective dual
varieties of the projectivized nilpotent
orbit closures in the isotropy modules of
symmetric spaces: our classi\-fi\-cation
yields its solution for all
$(-1)$-distingui\-shed orbits. To the best
of our know\-ledge, at this writing this
problem is largely open; even finding the
dimensions of these projective dual
varieties would be interesting (for the
minimal nilpotent orbit closures in the
adjoint modules, a solution was found in
\cite{KM}; cf.\;\cite{Sn} for a short
conceptual proof).

\section*{2. Main reductions
\label{formulation}}

%%In this section, after recalling a
%%necessary information from \cite{P2}, we
%%formulate some of our main results. In
%%Sections 3 and 4 we use them for obtaining
%%explicit lists and tables classifying
%%$(-1)$-distinguished elements.

Let $\g$ be a semisimple complex Lie
algebra, let $G$ be the adjoint group of
$\g$, and let $\theta \in {\rm Aut}\,\g$
be an element of order 2. We set
\begin{equation}
\ka:=\{x\in \g\mid \theta(x)=x\},\quad
\pe:=\{x\in \g\mid \theta(x)=-x\}.
\label{kp}
\end{equation}
Then $\ka$ and $\pe\neq 0$, the subalgebra
$\ka$ of~$\g$ is reductive, and
$\g=\ka\oplus \pe$ is a $\mathbb
Z_2$-grading of the Lie algebra $\g$,
cf.,\;e.g.,\;\cite{OV}. Denote by $G$ the
adjoint group of~$\g$. Let $K$ be the
connected reductive subgroup of $G$ with
the Lie algebra $\ka$; it is the adjoint
group of $\ka$.

Consider the adjoint $G$-module $\g$. Its
null-cone $\cNg$ is the cone of all
nilpotent elements of $\g$; it contains
only finitely many $G$-orbits, \cite{D},
\cite{K1}. The space $\pe$ is $K$-stable,
and we have
$$\cNp=\cNg\cap \pe.$$
By \cite{KR}, there are only finitely many
$K$-orbits in $\cNp$ and
\begin{equation}\label{2} 2\dim^{}_{\C} K\cdot x
=\dim^{}_{\C} G\cdot x \text{ \rm for any
} x\in \cNp.
\end{equation}

The Killing form $B_\g$ of $\g$
(resp.,\;\hskip -.4mm\;its restriction
$B_{\g}|{}_{\pe}$ to $\pe$) is symmetric
nondegenerate and $G$-invariant
(resp.,\;\hskip -.4mm\;$K$-invariant). We
identify the linear spaces $\g$ and $\g^*$
(resp.,\;\hskip -.4mm\;$\pe$ and~$\pe^*$)
by means of $B_\g$ (resp.,\;\hskip
-.4mm\;$B_{\g}|{}_{\pe}$).

Recall from \cite{BaCa}, \cite{Ca} that a
nilpotent element $x\!\in\! \cNg$ or its
$G$-orbit $G\cdot x$ is cal\-led {\it
distinguished } if the cent\-ra\-lizer
${\z}_{\g}(x)$ of $x$ in $\g$ contains no
nonzero semisimple elements. According to
\cite{P2}, distinguished orbits admit the
following geometric characterization:
\begin{theorem} {\rm (\cite[Theorem 2]{P2})}
\label{dist} Let $\cal O$ be a nonzero
$G$-orbit in $\cNg$ and $X:=\mathbf
P(\ov{\cal O})$. The following properties
are equivalent: \blist \item[\rm (i)]
$X=\check X$, \item[\rm (ii)] $\cal O$ is
distinguished.
\end{list}
\end{theorem}

Hence the classification of distinguished
elements obtained in \cite{BaCa} yields
the classification of projective self-dual
orbit closures in ${\mathbf P}(\cNg)$.
This, in turn, helps studying geometric
properties of these projective self-dual
projective varieties, \cite{P2}.

There is a counterpart of Theorem
\ref{dist} for the action of $K$ on $\pe$.
%%Namely, for any element $x\in \g$, denote
%%by ${\z}_{\g}(x)$ the centralizer of $x$
%%in $\g$.
Namely, if $x\in \pe$, then ${\z}_{\g}(x)$
is $\theta$-stable, hence it is a graded
Lie subalgebra of $\g$,
$$
{\z}_{\g}(x)={\z}_{\ka}(x)\oplus
{\z}_{\pe}(x),\ \mbox{ where } \
{\z}_{\ka}(x):={\z}_{\g}(x)\cap\ka, \ \
{\z}_{\pe}(x):={\z}_{\g}(x)\cap \pe.
$$
\begin{definition} {\rm (\cite{P2})}
\label{-1} {\rm An element $x\in \cNp$ and
its $K$-orbit are called}
$(-1)$-dis\-tinguished  {\rm if
${\z}_{\pe}(x)$ contains no nonzero
semisimple elements.}
\end{definition}

%%\begin{remark}
\noindent{\it Remark.}  It is tempting to
use simpler terminology and replace
``$(-1)$-distinguished'' with merely
``distinguished''. However this might lead
to the ambiguity since, as it is shown
below, there are the cases
%%(see
%%tables in Section 4)
where $x$ is $(-1)$-distinguished with
respect to the action of $K$ on $\p$ but
$x$ is not distinguished with respect to
the action of $G$ on $\g$.

\medskip

Notice that since $\pe$ contains nonzero
semisimple elements (see below the proof
of Theorem \ref{criter}), every
$(-1)$-distinguished element is nonzero.

\begin{theorem}\label{-dist}
{\rm (\cite[Theorem 5]{P2})} Let $\cal O$
be a nonzero $K$-orbit in $\cNp$ and
$X:=\mathbf P(\ov{\cal O})$. The following
properties are equivalent: \blist
\item[\rm (i)] $X=\check X$. \item[\rm
(ii)] $\cal O$ is $(-1)$-distinguished.
\end{list}
\end{theorem}

By virtue of Theorem \ref{-dist}, our goal
is classifying $(-1)$-distinguished
orbits.

\begin{theorem}\label{levi} Let $\cal O$ be
a $K$-orbit in $\cNp$ and let $x$ be an
element of $\cal O$. The following
properties are equivalent: \blist
\item[\rm (i)] $\cal O$ is
$(-1)$-distinguished.
 \item[\rm (ii)] The reductive Levi factors
of ${\z}_{\ka}(x)$ and ${\z}_{\g}(x)$ have
the same dimension. \item[\rm (iii)] The
reductive Levi factors of ${\z}_{\ka}(x)$
and ${\z}_{\g}(x)$ are isomorphic.
\end{list}
\end{theorem}

If $\g=\mathfrak h\oplus\mathfrak h$ where
$\mathfrak h$ is a semisimple complex Lie
algebra, and $\theta((y, z))=(z, y)$, then
\begin{equation}\label{reduction}
\ka=\{(y, y)\mid y\in \mathfrak h\},\quad
\p= \{(y, -y)\mid y\in \mathfrak h\}.
\end{equation}
If $H$ is the adjoint group of $\mathfrak
h$, then \eqref{reduction} implies that
$K$ is isomorphic to $H$ and the
$K$-module $\p$ is isomorphic to the
adjoint $H$-module $\mathfrak h$.
Therefore, for $y\in {\cal N}(\mathfrak
h)$,
  the variety ${\mathbf
P}(\ov{K\cdot y})$ is projective self-dual
if and only if the variety ${\mathbf
P}(\ov{H\cdot x})$ for $x:=(y, -y)\in
\cNp$ is projective self-dual. This is
consistent with Theorems \ref{dist} and
\ref{-dist}: the definitions clearly imply
that $x$ is $(-1)$-distinguished if and
only $y$ is distinguished. So Theorem
\ref{dist} follows from Theorem
\ref{-dist}.

In \cite{N}, it was made an attempt to
develop an analogue of the Bala--Carter
theory for nilpotent orbits in real
semisimple Lie algebras. The
Kostant--Sekiguchi bijection (recalled
below in this section) reduces this to
finding an analogue of the Bala--Carter
theory for $K$-orbits in $\cNp$. Given the
above discussion,
%%Theorems
%%\ref{dist} and \ref{-dist},
we believe that if such a natural analogue
exists, $(-1)$-distinguished elements
should play a key role in it, analogous to
that of distinguished elements in the
original Bala--Carter theory (in \cite{N},
the so called noticed elements, different
from $(-1)$-distinguished ones, play a
central role).

Theorems \ref{-dist}, \ref{levi} and their
proofs are valid over any algebraically
closed ground field of characteristic
zero. They give a key to explicit
classifying $(-1)$-distinguished elements
since there is a method (see the end of
this section) for finding reductive Levi
factors of ${\z}_{\ka}(x)$ and
${\z}_{\g}(x)$.

However, over the complex numbers, there
is another nice characterization of
$(-1)$-dis\-tinguished orbits. It is given
below in Theorem \ref{kost-sek}. For
classical $\g$, our approach to explicit
classifying $(-1)$-dis\-tinguished orbits
is based on this result. The very
formulation of this characterization uses
complex topology. However remark that, by
the Lefschetz's prin\-ciple, the final
explicit classification of
$(-1)$-dis\-tinguished orbits for
classical $\g$ given by combining Theorem
\ref{kost-sek} with Theorems
\ref{sl-compa}--\ref{skew-HermH-compa} is
valid over any algebraically closed ground
field of characteristic zero.

 Namely, according to the classical theory
 (cf.,\;e.g.,\;\cite{OV}), there is a $\theta$-stable real form $\gR$ of
$\g$ such that
\begin{equation}
\gR=\kR\oplus\pR,\ \mbox{ where }\
\kR:=\gR\cap\ka, \ \pR:=\gR\cap \pe,
\label{decom}
\end{equation}
is a Cartan decomposition of $\gR$. The
semisimple real algebra $\gR$ is
noncompact since $\pe\neq 0$. Assigning
$\gR$ to $\theta$ induces well defined
bijection from the set of all conjugacy
classes of elements of order $2$ in ${\rm
Aut}\,\g$ to the set of conjugacy classes
of noncompact real forms of~$\g$.  If $x$
is an element of $\gR$, we put
$$
{\z}_{\gR}(x):=\gR\cap {\z}_{\g}(x);
$$
${\z}_{\gR}(x)$ is the centralizer of $x$
in $\gR$. The identity component of the
Lie group of real points of $G$ is the
adjoint group ${\rm Ad}(\g^{}_{\R})$
of~$\gR$. We set
$$
\cNgr:=\cNg\cap \gR.
$$
\begin{definition} {\rm An element $x\in
\cNgr$ is called} compact {\rm if the
reductive Levi factor of its centralizer
${\z}_{\gR}(x)$ is a compact Lie algebra.}
\end{definition}

Recall that there is a bijection between
the sets of nonzero $K$-orbits in $\cNp$
and nonzero ${\rm Ad}(\g^{}_{\R})$-orbits
in $\cNgr$, cf.\,\cite{CM}, \cite{M}.
%%, called the
%%Kostant--Sekiguchi bijection.
Namely, let $\sigma$ be the complex
conjugation of $\g$ defined by $\gR$,
viz.,
$$\sigma(a+ib)=a-ib, \  a, b\in \gR.$$
Let $\{e, h, f\}$ be an {\it
$\mathfrak{sl}_2$-triple} in $\g$, i.e.,
an ordered triple of elements of $\g$
spanning a three-dimensional simple
subalgebra of $\g$ and satisfying the
bracket relations
\begin{equation}\label{bracket}
[h, e]=2e,\ [h, f]=-2f, \ [e, f]=h.
\end{equation}
It is called a {\it complex Cayley triple}
if $e, f\in \pe$ (hence $h\in \ka$) and
$\sigma (e)=-f$. Given a complex Cayley
triple $\{e, h, f\}$, set
\begin{equation}e':=
i(-h+e+f)/2, \quad h':=e-f,\quad f':=-
i(h+e+f)/2. \label{cayley1}
\end{equation}
Then $\{e', h', f'\}$ is an
$\mathfrak{sl}_2$-triple in $\gR$ such
that $\theta(e')=f'$ (and hence
$\theta(f')=e', \theta(h')=-h'$).
%%One can
%%show that $$G\!\cdot\!e=G\!\cdot e'.$$
The $\mathfrak{sl}_2$-triples in $\gR$
satisfying the last property are called
the {\it real Cayley triples}. Given
$\{e', h', f'\}$, the following formulas
restore $\{e, h, f\}$:
\begin{equation}e:=
(h'+if'-ie')/2, \quad h:=i(e'+f'),\quad
f:= (-h'+if'-ie')/2. \label{cayley2}
\end{equation}
The map $\{e, h, f\}\mapsto \{e', h',
f'\}$ is a bijection from the set of
complex to the set of real Cayley triples.
The triple $\{e, h, f\}$ is called the
{\it Cayley transform} of $\{e', h',
f'\}$.

Now let ${\cal O}$ be a nonzero $K$-orbit
in ${\cal N}(\pe)$. Then, according to
\cite{KR}, there exists a complex Cayley
triple $\{e, h, f\}$ in $\g$ such that
$e\in {\cal O}$. Let $\{e', h', f'\}$ be
the real Cayley triple in $\gR$
 such that $\{e, h, f\}$ is its Cayley
transform. Let ${\cal O}'=
%%$
%%be the
%%$
\GRc\!\cdot\!e'$.
%%-orbit
%%of $e'$.
Then the map assigning ${\cal O}'$ to
$\cal O$ is well defined and establishes a
bijection, called the {\it
Kostant--Sekiguchi bijection}, between the
set of nonzero $K$-orbits in ${\cal
N}(\pe)$ and the set of nonzero
$\GRc$-orbits in $\cNgr$. We say that the
orbits $\cal O$ and $\cal O'$ {\it
correspond} one another via the
Kostant--Sekiguchi bijection. One can
show, \cite{KR}, cf.\,\cite{CM}, that
\begin{gather}\label{dimKS}
\begin{gathered}
G\!\cdot\!e'=G\!\cdot e,\\
\dim_{\R}
\GRc\!\cdot\!e'=\dim_{\C}G\!\cdot\! e'.
\end{gathered}
\end{gather}

\begin{theorem}\label{kost-sek}
Let $\cal O$ be a nonzero $K$-orbit in
$\cNp$ and let $x$ be an element of the
$\GRc$-orbit in $\cNgr$ corresponding to
$\cal O$ via the Kostant--Sekiguchi
bijection. The following pro\-perties are
equivalent: \blist \item[\rm (i)] $\cal O$
is $(-1)$-distinguished. \item[\rm (ii)]
$x$ is compact.
\end{list}
\end{theorem}

Next result, Theorem \ref{decomp}, reduces
studying projective self-dual varieties
associated with symmetric spaces to the
case of simple Lie algebras $\g$. Namely,
since the Lie algebra $\g$ is semisimple,
it is the direct sum of all its simple
ideals. As this set of ideals is
$\theta$-stable, $\g$ is the direct sum of
its $\theta$-stable ideals,
%%~$\g_i$,
$ \g=\g_1\oplus\ldots\oplus\g_d, $ where
each $\g_l$ is either \blist \item[\rm
(a)] simple or \item[\rm (b)] direct sum
of two isomorphic simple ideals permuted
by $\theta$.\end{list}

Let $G_l$, $\ka_l$, $\p_l$, $K_l$ and
${\cal N}(\p_l)$ have the same meaning for
$\g_l$ with respect to
$\theta_l:=\theta|{}_{\g_l}$ as
resp.,\;\hskip -.4mm\;$G$, $\ka$, $\p$,
$K$ and ${\cal N}(\p)$ have for $\g$ with
respect to $\theta$. Then
$G=G_1\times\ldots\times G_d$,
$\ka=\ka_1\oplus\ldots\oplus\ka_d$,
$\p=\p_1\oplus\ldots\oplus\p_d$,
$K=K_1\times\ldots\times K_d$ and ${\cal
N}(\p)= {\cal N}(\p_1)\times\ldots\times
{\cal N}(\p_d)$.

In Theorem \ref{decomp}, we use the
following notation. Let $X_1,\ldots, X_s$
be the closed sub\-varieties of a
projective space $\mathbf P$, and let
$s\geqslant 2$. Consider the variety
\begin{equation}
\overline{\{(x_1,\ldots, x_s, y)\in
X_1\times\ldots\times X_s\times \mathbf
P\mid \dim \langle x_1,\ldots, x_s\rangle=
s-1, \ y\in \langle x_1,\ldots,
x_s\rangle\}}, \label{join}
\end{equation}
where bar denotes Zariski closure in
$X_1\times \ldots\times X_s\times \mathbf
P$ and $\langle S\rangle$ denotes the
linear span of $S$ in $\mathbf P$.
Consider the projection of variety
$\eqref{join}$ to $\mathbf P$. Its image
is denoted~by
\begin{equation}
{\rm Join}(X_1,\ldots, X_s) \label{join1}
\end{equation}
and called the {\it join} of $X_1,\ldots,
X_s$. If $s=2$, then \eqref{join1} is the
usual join of $X_1$ and $X_2$,
cf.\,\cite{Ha}. If $s> 2$, then
\eqref{join1} is the usual join of ${\rm
Join}(X_1,\ldots, X_{s-1})$ and
 $X_s$.

\begin{theorem}\label{decomp}
Let  $x=x_1+\ldots+x_d$ where $x_l\in
{\cal N}(\p_l)$, $l=1,\ldots, d$. Consider
in $\mathbf P(\g)$ the projective
subvarieties $X:=\mathbf P(\ov{K\cdot x})$
and $X_l:=\mathbf P(\ov{K_l\cdot x_l})$.
Then
$$ X={\rm Join}(X_1,\ldots, X_d)$$ and the
following properties are equivalent:
\blist \item[\rm (i)] $X=\check X$.
\item[\rm (ii)] $X_l=\check X_l$ for all
$l$.
 \end{list}
\end{theorem}

\noindent{\it Remark.} This is a specific
geometric property of the considered
varieties. In general setting, projective
self-duality of  ${\rm Join}(Z_1,\ldots,
Z_m)$ is not equivalent to projective
self-duality of all $Z_1,\ldots, Z_m$.

\medskip

If $\g_l$ is of type (b), viz., $\g_l$ is
isomorphic to the direct sum of algebras
$\mathfrak s\oplus \mathfrak s$, where
$\mathfrak s$ is a simple algebra, and
$\theta_l$ acts by $\theta_l((x, y))=(y,
x)$, then, as explained above (see
\eqref{reduction}), classifying
$(-1)$-distinguished orbits in ${\cal
N}(\p_l)$ amounts to classifying
distinguished orbits in ${\cal
N}(\mathfrak s)$. Since the last
classification is known, \cite{BaCa},
\cite{Ca}, this and Theorem \ref{decomp}
reduce  classifying $(-1)$-distinguished
orbits in $\g$ to the case where $\g$ is a
{\it simple} Lie algebra.

In this case, for explicit classifying
$(-1)$-distinguished orbits, one can apply
Theorem~\ref{levi} since there are
algorithms, \cite{Ka}, \cite{Vi},
yielding, in principle, a classification
of $K$-orbits in ${\cal N}(\pe)$ and a
description of the reductive Levi factors
of $\z_{\ka}(x)$ and $\z_{\g}(x)$ for
$x\in{\cal N}(\pe)$. Moreover, there is a
computer program, \cite{L}, implementing
these algorithms
 and yielding, for concrete pairs
 $(\g, \theta)$,
 the explicit
representatives $x$ of $K$-orbits in
$\cNp$ and the reductive Levi factors of
$\z_{\ka}(x)$ and $\z_{\g}(x)$.
%%In principle,
There is also another way for finding
representatives and dimensions of
$K$-orbits in $\cNp$: they may be obtained
using
%%can be also
%%obtained by means of the
the general algorithm of finding Hesselink
strata for linear reductive group actions
given in \cite{P3}, since in our case
these strata coincide with $K$-orbits, see
\cite[Proposition 4]{P3}.

Fortunately, for exceptional simple $\g$,
explicit finding of the aforementioned
classification and description already has
been  performed by D.~\v Z.~{--}~\hskip
-2.8mm Do\-ko\-vi\'c in \cite{D3},
\cite{D4}, cf.\;\cite{CM}.  In these
papers, the answers are given in terms of
the so called {\it characteristics} in the
sense of Dynkin (in Section 5, we recall
this classical approach, \cite{D},
\cite{K1}, \cite{KR}, cf.\;\cite{CM} and
\cite{Vi}). Combining either of Theorems
\ref{levi} and
%%or Theorem
\ref{kost-sek} with these results, we
obtain, for all exceptional simple $\g$
and all $\theta$, the explicit
classification of $(-1)$-distinguished
$K$-orbits in $\cNp$ in terms of their
characteristics (see Theorem
\ref{except}).

For classical $\g$, we use another
approach and obtain the classification of
$(-1)$-distin\-gui\-shed $K$-orbits in
$\cNp$ by means of the Kostant--Sekiguchi
bijection, Theorem \ref{kost-sek} and
elementary representation theory of
$\ssl_2$ (see Theorems
\ref{sl-compa}--\ref{skew-HermH-compa}).

%%\medskip

Theorems \ref{levi}, \ref{kost-sek},
\ref{decomp} are proved in Section 3. For
simple Lie algebras $\g$, the
classifications of $(-1)$-distingui\-shed
$K$-orbits in $\cNp$ are obtained in the
next two sections: in Section~4, we
consider the case of classical $\g$, and
in Section 5, that of exceptional simple
$\g$. In Section~6, we briefly discuss
some geometric properties of the
constructed self-du\-al projective
algebraic varieties.

\section*{3. Proofs of Theorems 4--6
}

%%\subsection*{Proof of Theorem \ref{levi}}
{\it Proof of Theorem \ref{levi}.} Theorem
\ref{levi} immediately follows from
Theorem \ref{criter} proved below.

 Let $k$ be a
 field of characteristic $0$.
 If $\mathfrak c$ is a finite dimensional
 algebraic Lie algebra over $k$, we denote
 by $\rd\,\mathfrak c$ (resp.,\;\hskip -.4mm\;$\ru\mathfrak
 c$) the radical (resp.,\;\hskip -.4mm\;the unipotent
 radical, i.e., the maximal ideal
 whose elements are
 nilpotent) of $\mathfrak c$. We call
 $\mathfrak c/\rd\,\mathfrak c$
 (resp.,\;\hskip -.4mm\;$\mathfrak c/\ru\mathfrak c$)
 the {\it
Levi
 factor} (resp.,\;\hskip -.4mm\;the {\it reductive
Levi factor})
 of $\mathfrak c$. According to \cite{Ch}, there
 is a semisimple (resp.,\;\hskip -.4mm\;reductive) subalgebra
 $\mathfrak l$ in
 $\mathfrak c$ such that $\mathfrak c$ is
 the se\-mi\-direct
 sum of $\mathfrak l$ and $\rd\,\mathfrak
 c$ (resp.,\;\hskip -.4mm\;$\ru\mathfrak c$). Every such
 $\mathfrak l$ is called {\it Levi
 subalgebra}
 (resp.,\;\hskip -.4mm\;{\it reductive Levi subalgebra}) of
 $\mathfrak c$.

\begin{theorem} \label{criter}
Assume that $k$ is algebraically closed
and let $\mathfrak a$ be a finite
dimensional algebraic Lie algebra over
$k$. Let
$\theta\in\operatorname{Aut}\mathfrak a$
be an element of order $2$. For any
$\theta$-stable linear subspace~$\mathfrak
l$ of $\mathfrak a$, set
$$\mathfrak l^{\pm}\!=\!\{x\!\in\!
\mathfrak l\mid \theta(x)\!=\!\pm x\}.$$
The following properties are
equivalent{\rm :} \blist \item[\rm(i)]
$\mathfrak a^-$ contains no nonzero
semisimple elements. \item[\rm(ii)]
$\mathfrak a^-=(\ru\mathfrak a)^-$.
\item[\rm(iii)] $\dim \mathfrak a^-=\dim
(\ru\mathfrak a)^-$. \item[\rm(iv)]
$\mathfrak a^-\subseteq \ru\mathfrak a$.
\item[\rm (v)] The reductive Levi factors
of ${\mathfrak a}^{+}$\!\! and $\mathfrak
a$ have the same dimension. \item[\rm
(vi)] The reductive Levi factors of
$\mathfrak a^+$\!\! and $\mathfrak a$ are
isomorphic. \item[\rm(vii)] The set of all
reductive Levi subalgebras of $\mathfrak
a^+$\!\! coincides with the set of all
$\theta$-stable reductive Levi subalgebras
of $\mathfrak a$.
\end{list}
\end{theorem}
\noindent {\it Proof.} It is known that
$\mathfrak a$ contains a $\theta$-stable
Levi subalgebra $\mathfrak s$
(e.\,g.,\;see\;\cite[Appen\-dix~9]{KN}),
and $\rd\,\mathfrak a$ contains a
$\theta$-stable maximal torus $\mathfrak
t$, see\;\cite{St}. Thus we have the
fol\-lowing $\theta$-stable direct sum
decompositions of vector spaces:
\begin{gather}
\begin{cases}
\begin{gathered}
{\mathfrak a}={\mathfrak
s}\oplus{\rd}\,{\mathfrak a}, \quad
{\rd}\,{\mathfrak a}={\mathfrak
t}\oplus {{\rm rad}_u} {\mathfrak a},\\[4pt]
\mathfrak s=\mathfrak s^+\oplus \mathfrak
s^-,\quad \mathfrak t=\mathfrak t^+\oplus
\mathfrak t^-, \quad \ru\mathfrak
a=(\ru\mathfrak
a)^+\oplus (\ru\mathfrak a)^-,\\[4pt]
\mathfrak a^\pm=\mathfrak
s^\pm\oplus\mathfrak
t^\pm\oplus(\ru\mathfrak a)^\pm.
\label{ld}
\end{gathered}
\end{cases}
\end{gather}

${\rm (i)}\Rightarrow{\rm (ii)}$: Assume
that (i) holds. Then, by \eqref{ld}, we
have $\mathfrak t^-=0$. If $\mathfrak
s^-\neq 0$, then $\theta|^{}_{\mathfrak
s}\in \operatorname{Aut}\,\mathfrak a$ is
an element of order 2. By \cite{Vu}, this
implies that $\mathfrak s^-$ contains a
nonzero $\theta$-stable algebraic torus.
By (11), this contradicts (i). Thus
$\mathfrak s^-=0$. Whence (ii)
by~\eqref{ld}.

${\rm (ii)}\Rightarrow{\rm (i)}$ and ${\rm
(iv)}\Rightarrow{\rm (i)}$: This is
because all elements of $\ru\mathfrak a$
are nilpotent.

${\rm (ii)}\Rightarrow{\rm (iii)}$ and
${\rm (ii)}\Rightarrow{\rm (iv)}$: This is
evident.

${\rm (iii)}\Rightarrow{\rm (ii)}$: This
follows from \eqref{ld}.

${\rm (ii)}\Rightarrow{\rm (v)}$: It
follows from \eqref{ld} that (ii) is
equivalent to
\begin{equation}
\mathfrak s\oplus\mathfrak
t\subseteq\mathfrak a^+.\label{ld1}
\end{equation}
Assume that \eqref{ld1} holds. Let $\pi:
\mathfrak a\rightarrow\mathfrak
a/\ru\mathfrak a$ be the natural
homomorphism. By \eqref{ld}, \eqref{ld1},
we have $\pi(\mathfrak a^+)=\mathfrak
a/\ru\mathfrak a$. Since the algebra
$\mathfrak a/\ru\mathfrak a$ is reductive,
this implies
\begin{gather}
\begin{gathered}
\dim \mathfrak a/\ru\mathfrak
a\leqslant\mbox{dimension of reductive
Levi
subalgebras of} \ \mathfrak a^+\\
=\dim \mathfrak a^+/\ru(\mathfrak a^+).
\end{gathered}
\label{ld2}
\end{gather} On the other hand, since
reductive Levi subal\-gebras of $\mathfrak
a$ and $\mathfrak a^+$ are precisely their
maximal reductive subal\-gebras,
cf.,\,e.g.,\,\cite{OV} , the inclusion
$\mathfrak a^+\subseteq \mathfrak a$
implies
\begin{gather}
\begin{gathered}
\dim \mathfrak a/\ru\mathfrak
a=\mbox{dimension of reductive Levi
subalgebras of }\ \mathfrak a\\
\geqslant\mbox{dimension of reductive Levi
subalgebras of }\ \mathfrak a^+.
\end{gathered}
\label{ld3}
\end{gather}
Now (v) follows from \eqref{ld2},
\eqref{ld3}.

${\rm (v)}\Rightarrow{\rm (ii)}$: Assume
that (v) holds. Let $\mathfrak l$ be a
reductive Levi subalgebra of $\mathfrak
a^+$. Since reductive Levi subalgebras of
$\mathfrak a$ are precisely its maximal
reductive subalgebras, it follows from (v)
that $\mathfrak l$ is a reductive Levi
subalgebra of $\mathfrak a $ as well. Let
$\mathfrak s$ and $\mathfrak t$ be
resp.\;\hskip -.4mm the derived subalgebra
and the center of $\mathfrak l$; so we
have $\mathfrak l=\mathfrak
s\oplus\mathfrak t$. Then $\mathfrak s$ is
a $\theta$-stable Levi subalgebra of
$\mathfrak a$, and $\mathfrak t$ is a
$\theta$-stable maximal torus of
$\rd\,\mathfrak a$. Hence, as above, (ii)
is equivalent to the inclusion
\eqref{ld1}. As the latter holds by the
very construction of $\mathfrak l$, we
conclude that (ii) holds as well.

${\rm (vii)}\Rightarrow{\rm (v)}$: This is
clear.

${\rm (i)}\Rightarrow{\rm (vii)}$: Let
$\mathfrak l$ be a $\theta$-stable
reductive Levi subalgebra of $\mathfrak
a$. If $\mathfrak l^-\neq 0$, then
$\theta^{}|_{\mathfrak l}\in {\rm Aut}\,
\mathfrak l$ is an element of order 2.
Hence, by \cite{Vu}, there is  a nonzero
$\theta$-stable algebraic torus in
$\mathfrak l^-$. This contradicts (i).
Whence $\mathfrak l^-=0$, i.e., $\mathfrak
l$ lies in~$\mathfrak a^+$.

${\rm (vi)}\Rightarrow{\rm (v)}$ and ${\rm
(vii)}\Rightarrow{\rm (vi)}$: This is
clear. \quad $\square$

\bigskip

%%\subsection*{Proof of Theorem \ref{kost-sek}}

\noindent{\it Proof of Theorem
\ref{kost-sek}.} According to \cite{KR},
there is a complex Cayley triple $\{e, h,
f\}$ in $\g$ such that ${\cal O}=K\cdot
e$. Let $\{e', h', f'\}$ be the real
Cayley triple in $\gR$ whose Cayley
transform is $\{e, h, f\}$. Since ${\cal
O}'=\GRc\cdot e'$, we may, and will,
assume that $x=e'$.

Let $\mathfrak s$ (resp.,\;\hskip
-.4mm\;$\mathfrak s^{}_{\R}$) be the
simple three-dimensional subalgebra of
$\g$ (resp.,\;\hskip -.4mm $\gR$) spanned
by $\{e, h, f\}$ (resp.,\;\hskip
-.4mm\;$\{e', h', f'\}$),
\begin{equation}
\mathfrak s=\C e+\C h+\C f,\quad \mathfrak
s_{\R}=\R e' +\R h' + \R f'.
\label{triples}
\end{equation}
Denote  by $\z_{\g}(\mathfrak s)$
(resp.,\;\hskip -.4mm\;$\z_{\gR}(\mathfrak
s^{}_{\R})$) the centralizer of $\mathfrak
s$ (resp.,\;\hskip -.4mm\;$\mathfrak
s^{}_{\R}$) in $\g$ (resp.,\;\hskip
-.4mm\;$\gR$).

It follows from \eqref{cayley1},
\eqref{cayley2}, \eqref{triples} that
$\mathfrak s^{}_{\R}$ is a real form of
$\mathfrak s$. Since $\gR$ is a real form
of $\g$, this yields that
$\z^{}_{\gR}(\mathfrak s^{}_{\R})$ is a
real form of $\z_{\g}(\mathfrak s)$. The
definitions of complex and real Cayley
triples imply that $\mathfrak s$ and
$\mathfrak s^{}_{\R}$ are $\theta$-stable
subalgebras. Therefore
$\z^{}_{\g}(\mathfrak s)$ and
$\z_{\gR}(\mathfrak s^{}_{\R})$ are
$\theta$-stable as well. Whence
\begin{gather}
\z^{}_{\g}(\mathfrak
s)=(\z^{}_{\g}(\mathfrak s)\cap \ka)\oplus
(\z^{}_{\g}(\mathfrak s)\cap \p), \quad
\z^{}_{\gR}(\mathfrak s^{}_{\R})=
(\z^{}_{\gR}(\mathfrak s^{}_{\R})\cap
\kR)\oplus (\z^{}_{\gR}(\mathfrak
s^{}_{\R})\cap
\pR),\label{c1}\\[2pt]
\z^{}_{\g}(\mathfrak s)\cap \ka=
\C(\z^{}_{\gR}(\mathfrak s^{}_{\R})\cap
\kR), \quad \z^{}_{\g}(\mathfrak s)\cap
\p=\C(\z^{}_{\gR}(\mathfrak s^{}_{\R})\cap
\pR).
 \label{c2}
\end{gather}

Now take into account that
$\z_{\g}(\mathfrak s)$ (resp.,\;\hskip
-.4mm\;$\z_{\gR}(\mathfrak s^{}_{\R})$) is
the reductive Levi subalgebra of
$\z_{\g}(e)$ (resp.,\;\hskip
-.4mm\;$\z_{\gR}(e')$): regarding
$\z_{\g}(\mathfrak s)$, see,
e.g.,\;\cite{SpSt}; for
$\z_{\gR}(\mathfrak s^{}_{\R})$, the
arguments are the same. Hence, by Theorem
\ref{criter}, the orbit $\cal O$ is
$(-1)$-distinguished if and only if
$\z_{\g}(\mathfrak s)\cap \p=0$. By
\eqref{c2}, the latter condition is
equivalent to $\z_{\gR}(\mathfrak
s^{}_{\R})\cap \pR=0$. On the other hand,
since \eqref{decom} is the Cartan
decomposition of $\gR$, it follows from
\eqref{c1} that $\z_{\gR}(\mathfrak
s^{}_{\R})\cap \pR=0$ if and only if the
Lie algebra $\z_{\gR}(\mathfrak
s^{}_{\R})$ is compact. This completes the
proof of Theorem \ref{kost-sek}.\quad
$\square $

%%\nopagebreak

\bigskip

%%\nopagebreak

%%\subsection*{Proof of Theorem \ref{decomp}}

\noindent {\it Proof of Theorem
\ref{decomp}.} Clearly the following
proposition implies the first statement.

\begin{proposition}
{\rm (\cite[Proposition 2]{P2})} Let $H_t$
be an algebraic group and let $L_t$ be a
finite dimensional algebraic $H_t$-module,
$t=1, 2$. Let $v_t\in L_s$ be a nonzero
vector such that the orbit $H_t\cdot v_t$
is stable with respect to scalar
multiplications. Put $H:=H_1\times H_2$,
$L:=L_1\oplus L_2$. Identify $L_t$ with
the linear subspace of $L$ and set
$v:=v_1+v_2$. Then
$$
\mathbf P(\overline{H\cdot v})={\rm
Join}(\mathbf P(\overline{H_1\cdot v_1}),
\mathbf P(\overline{H_2\cdot v_2})).
$$
\end{proposition}

Regarding the second statement, notice
that the centralizers ${\z}_{\g_l}(x_l)$
and ${\z}_{\g}(x)$ are resp.\;\hskip -.4mm
$\theta_l$- and $\theta$-stable and
${\z}_{\pe}(x)={\z}_{{\pe}_1}(x_1)\oplus\ldots\oplus
{\z}_{{\pe}_d}(x_d)$. Now the claim
follows from this decomposition,
Definition 1 and Theorem \ref{-dist}.\quad
$\square$

\section*{4. Classification for classical
\boldmath$\g$ }

If $\g$ is classical, our approach to
classifying $(-1)$-distinguished nilpotent
$K$-orbits in $\pe$ is based on the
Kostant--Sekiguchi bijection and Theorem
\ref{kost-sek}. Namely, for every
noncompact real form $\g^{}_{\R}$ of $\g$,
we will find in $\cNgr$
%%in $\g^{}_\R$
all
%%the
%%nilpotent
${\rm Ad}(\g^{}_{\R})$-orbits $\cal O$
%% of the adjoint group ${\rm Ad}(\g^{}_{\R})$
%%of $\g^{}_{\R}$
%%in $\cNgr$
 whose elements are compact.
%%such that the reductive Levi
%%factor of $\z^{}_{\g^{}_{\R}}(x)$ for
%%$x\in \cal O$ is compact.
If this is done, then for the
complexification
%%$\theta\in {\rm Aut}\,\g$
of a Cartan decomposition of $\g^{}_{\R}$,
the $K$-orbits in $\cNp$ corresponding to
all such $\cal O$ via the
Kostant--Sekiguchi bijection are precisely
all $(-1)$-distinguished orbits.

First we recall the classification of real
forms $\gR$ of classical complex Lie
algebras~$\g$. Let $D$ be either real
numbers $\R$, complex numbers $\C$, or
quaternions $\Ha$. If $a\in D$, we denote
by $\overline{a}$ the conjugate of $a$.
Let $V$ be a finite dimensional vector
space over $D$ (left $D$-module). Denote
by ${\GL}_D(V)$ the Lie group of all
$D$-linear automorphisms of $V$. Its
 Lie algebra $\gl_D(V)$ is
identified with the real Lie algebra of
all $D$-linear endomorphisms of $V$. The
derived algebra of $\gl_D(V)$ is denoted
by $\ssl_D(V)$.

Let a map $\Phi:\,V\times V\to D$ be
$D$-linear with respect to the first
argument. It is called \blist \item[---]
{\it symmetric bilinear form on} $V$ if
$\Phi(x,y)=\Phi(y,x)$ for all $x, y\in V$,
\item[---] {\it skew-symmetric bilinear
form on} $V$ if $\Phi(x,y)=-\Phi(y,x)$ for
all $x, y\in V$, \item[---] {\it Hermitian
form on} $V$ if $\Phi(x,y)=\ov{\Phi(y,x)}$
for all $x, y\in V$, \item[---] {\it
skew-Hermitian form on} $V$ if
$\Phi(x,y)=-\ov{\Phi(y,x)}$ for all $x,
y\in V$.
\end{list}

 By ``form on $V$" we always
mean form of one of these four types. A
form $\Phi$ on $V$ is called {\it
nondegenerate} if, for any $x\in V$, there
is $y\in V$ such that $\Phi(x,y)\ne 0$. If
$\Phi_s$, $s=1, 2$, is a form on a finite
dimensional vector space $V_s$ over $D$,
then $\Phi_1$ and $\Phi_2$ are called {\it
equivalent forms} if there exists an
isometry $\psi: V_1\rightarrow V_2$ with
respect to $\Phi_1$ and $\Phi_2$, i.e., an
isomorphism of vector spaces over $D$ such
that $\Phi_2(\psi(x), \psi(y))=\Phi_1(x,
y)$ for all $x, y\in V_1$.

The automorphism group of a form $\Phi$ on
$V$,
$${\GL}_D^{\Phi}(V):=\{g\in{\GL}_D(V)
\mid \Phi(g\cdot x, g\cdot y)=\Phi(x, y)\
\text{for all}\ x,y\in V\},$$ is a Lie
subgroup of ${\GL}_D(V)$ whose Lie algebra
is
$$\gl_D^{\Phi}(V):=\{A\in\gl_D(V)\mid
\Phi(A\cdot x,y)+\Phi(x,A\cdot y)=0\
\text{for all}\ x,y\in V\}.$$ For $\Phi=0$
(zero form), we have
${\GL}_D^{\Phi}(V)={\GL}_D(V)$ and
$\gl_D^{\Phi}(V)=\gl_D(V)$.

Let $n:=\dim _D V$ and let
$(e):=\{e_1,\ldots,e_n\}$ be a basis of
$V$ over $D$. Identifying $D$-linear
endomorphisms of $V$ with their matrices
with respect to $(e)$, we identify
$\gl_D(V)$, $\ssl_D(V)$ and
$\gl_D^{\Phi}(V)$ with the corresponding
matrix real Lie algebras. For $D=\R$,
$\Ha$, the real Lie algebras $\gl_D(V)$
and $\ssl_D(V)$ are denoted resp.\;\hskip
-.4mm by $\gl_n(D)$ and $\ssl_n(D)$, and
for $D=\C$, by $\gl_n(\C)_{\R}$ and
$\ssl_n(\C)_{\R}$. The last algebras are
endowed with the natural structures of
complex Lie algebras that are denoted
resp.\;\hskip -.4mm by $\gl_n(\C)$ and
$\ssl_n(\C)$. If $D=\Ha$, then $V$ is a
$2n$-dimensional vector space over the
subfield $\C$ of $\Ha$, and elements of
$\gl_{\Ha}(V)$ are its $\C$-linear
endomorphisms. Identifying them with their
matrices with respect to the basis
$e_1,\ldots, e_n, je_1,\ldots, je_n$, we
identify $\gl_n(\Ha)$ (resp.,\;\hskip
-.4mm\;$\ssl_n(\Ha)$) with the
corresponding Lie subalgebra of
$\gl_{2n}(\C)_{\R}$ (resp.,\;\hskip
-.4mm\;$\ssl_{2n}(\C)_{\R}$).

If $\Phi$ is a form on $V$, then the
matrix $(\Phi_{st})$ of $\Phi$ with
respect to $(e)$ defines $\Phi$ by
%%the equality
\begin{gather*}
\textstyle
 \Phi\bigl(\sum_{s=1}^n x_se_s,\sum_{t=1}^n
y_te_t\bigr)= \begin{cases} \textstyle
\sum_{s,t=1}^n x_s\Phi_{st}y_t \ \text{ if
$\Phi$ is symmetric or skew-symmetric},\\
 \textstyle \sum_{s,t=1}^n
x_s\Phi_{st}\ov{y_t}\ \text{ if $\Phi$ is
Hermitian or skew-Hermitian.}
\end{cases}
\end{gather*}

  The map $\Phi\mapsto (\Phi_{st})$
is a bijection between all symmetric
(resp.,\;\hskip -.4mm\;skew-symmetric,
Hermitian, skew-Hermitian) forms on $V$
and all symmetric (resp.,\;\hskip
-.4mm\;skew-symmet\-ric, Hermitian,
skew-Hermitian) $n\times n$-matrices over
$D$. $\Phi$ is nondegenerate if and only
if $(\Phi_{st})$ is nonsingular.

We denote by ${\rm I}_d$ the unit matrix
of size $d$ and put ${\rm I}_{p,q}:={\rm
diag} (1,\ldots,1, -1,\ldots,-1)$ where
$p$ (resp.,\;\hskip -.4mm\;$q$) is the
number of 1's (resp.,\;\hskip
-.4mm\;$-1$'s).

\smallskip

Now we summarize the results from linear
algebra about classification of forms and
fix some notation and terminology. In the
sequel, $\mathbb N$ denotes the set of all
nonnegative integers.

\vskip 2.5mm

{\it Symmetric bilinear forms}

\smallskip

$D=\R$. There are exactly $n+1$
equivalence classes of nondegenerate
symmetric bilinear forms on $V$. They are
represented by the forms $\Phi$ with
$(\Phi_{st})={\rm I}_{p,q}$, $p+q=n$,
$p=0,\dots, n$. If $(\Phi_{st})={\rm
I}_{p,q}$, then ${\rm sgn}\,\Phi:=p\!-\!q$
is called the {\it signature of} $\Phi$,
and the real Lie algebra $\gl_D^{\Phi}(V)$
is denoted by~$\so_{p,q}$.

$D=\C$. There is exactly one equivalence
class of nondegenerate symmetric bilinear
forms on $V$. It is represented by the
form $\Phi$ with $(\Phi_{st})={\rm I}_n$.
The corresponding real Lie algebra
$\gl_D^{\Phi}(V)$ is denoted by
$\so_n(\C)_{\R}$. It has a natural
structure of complex Lie algebra that is
denoted by $\so_n(\C)$.

$D=\Ha$. There are no nonzero symmetric
bilinear forms on $V$.

\medskip {\it Skew-symmetric bilinear forms}

\smallskip

$D=\R$ and $\C$. Nondegenerate
skew-symmetric bilinear forms on $V$ exist
if and only if $n$ is even. In this case,
there is exactly one equivalence class of
such forms. It is represented by the form
$\Phi$ with
$(\Phi_{st})=\Bigl[\begin{smallmatrix} 0&{\rm I}_{n/2}\\
-{\rm I}_{n/2}&0\end{smallmatrix}\Bigr]$.
Corresponding real Lie algebra
$\gl_D^{\Phi}(V)$ is denoted by
${\spl}_n(\R)$ for $D=\R$, and by
$\spl_n(\C)_{\R}$ for $D=\C$. The algebra
$\spl_n(\C)_{\R}$ has a natural
struc\-ture of complex Lie algebra that is
denoted by $\spl_n(\C)$.

$D=\Ha$. There are no nonzero
skew-symmetric bilinear forms on $V$.

\medskip {\it Hermitian forms}

\smallskip

$D=\R$. Hermitian forms on $V$ coincide
with symmetric bilinear forms.

$D=\C$ and $\Ha$. There are exactly $n+1$
equivalence classes of nondegenerate
Hermitian forms on $V$. They are
represented by the forms $\Phi$ with
$(\Phi_{st})={\rm I}_{p,q}$, $p+q=n$,
$p=0,\ldots, n$. If $(\Phi_{st})={\rm
I}_{p,q}$, then ${\rm sgn}\,\Phi:=p\!-\!q$
is called the {\it signature} of $\Phi$,
and the real Lie algebra $\gl_D^{\Phi}(V)$
is denoted by $\su_{p,q}$ for $D=\C$, and
by $\spl_{p,q}$ for $D=\Ha$.

\medskip {\it Skew-Hermitian forms}

\smallskip

$D=\R$. Skew-Hermitian forms on $V$
coincide with skew-symmetric forms.

$D=\C$. The map $\Phi\mapsto i\Phi$ is a
bijection between all nondegenerate
Hermitian forms on $V$ and all
nondegenerate skew-Hermitian forms on $V$.

$D=\Ha$. There is exactly one equivalence
class of nondegenerate skew-Hermitian
forms on $V$. It is represented by the
form $\Phi$ with $(\Phi_{st})=j{\rm I}_n$.
The corresponding real Lie algebra
$\gl_D^{\Phi}(V)$ is denoted by $\un
^*_n(\Ha)$.

\vskip 2.5mm

Some of the above-defined Lie algebras are
isomorphic to each other. Obviously
$\su_{p,q}=\su_{q,p}$,
$\so_{p,q}=\so_{q,p}$,
$\spl_{p,q}=\spl_{q,p}$ and we also have
\begin{gather}\label{iso}
\begin{cases}
\begin{gathered}
\so_2=\un_1^*(\Ha),\
\so_3(\C)\simeq\ssl_2(\C)=\spl_2(\C),\
\so_3\simeq \su_2=
\spl_{1,0}=\ssl_1(\Ha),\\[1pt]
\so_{1,2}\simeq
\su_{1,1}\simeq\ssl_2(\R)=\spl_2(\R),\
\so_4(\C)\simeq\ssl_2(\C)\oplus\ssl_2(\C),\
\so_4\simeq \su_2\oplus\su_2,\\[1pt]
\so_{1,3}\simeq\ssl_2(\C)_{\R},
\so_{2,2}\simeq\ssl_2(\R)\oplus\ssl_2(\R),
\un_2^*(\Ha)\simeq\su_2\oplus\ssl_2(\R),
\so_5(\C)\simeq\spl_4(\C),\\[1pt]
\so_5\simeq\spl_2,
\so_{1,4}\simeq\spl_{1,1},\
\so_{2,3}\simeq\spl_4(\R),\
\so_6(\C)\simeq\ssl_4(\C),\
\so_6\simeq\su_4,\\[1pt]
\so_{1,5}\simeq \ssl_2(\Ha),\
\so_{2,4}\simeq \su_{2,2},\
\so_{3,3}\simeq \ssl_4(\R),\
\un_3(\Ha)\simeq\su_{1,3},\ \un
_4(\Ha)\simeq\so_{2,6}.
\end{gathered}
\end{cases}
\end{gather}

By definition, {\it classical complex Lie
algebras} $\g$ are the algebras
$\ssl_n(\C)$, $\so_n(\C)$ and
$\spl_n(\C)$. According to E.\;Car\-tan's
classification (cf.,\;e.g.,\;\cite{OV}),
up to isomorphism, all their real forms
$\gR$ are listed in Table~1.

\

\vskip -11mm

\

\begin{center}{TABLE 1}\\[6pt]
{\fontsize{9pt}{5mm}\selectfont
\begin{tabular}{c|c|c}
\ $\hskip .2mm \g\hskip 1.5mm $ & $\gR$ &
compact $\gR$
\\
\hline \hline
&\\[-11pt]
$\ssl_n(\C),\ n\geqslant 2$&
$\begin{matrix}
\ssl_n(\R),\\
\ssl_l(\Ha), \ n=2l,\\
\su_{n-q,q}, \ q=0, 1,\ldots, [n/2]
\end{matrix}$ & $\su_{n}:=\su^{}_{n, 0}$
\\[-14pt]
&&
\\[3pt]
\hline
&\\[-12pt]
$\so_n(\C),\ n=3$ or $n\geqslant 5$&
$\begin{matrix} \so_{n-q, q}, \ q=0,
1,\ldots,
[n/2],\\
\mathfrak u^*_{l}(\Ha),\ n=2l
\end{matrix}$&$\so_{n}:=\so^{}_{n, 0}$
\\[-14pt]
&&
\\[3pt]
\hline
&\\[-12pt]
$\spl_n(\C),\ n=2l\geqslant 2$ &
$\begin{matrix}
\spl_n(\R),\\
\spl_{l-q, q}, \ q=0, 1,\ldots, [l/2]
\end{matrix}$&$\spl_{l}:=\spl^{}_{l, 0}$
\end{tabular}
}
\end{center}

\vskip 1mm

There are the isomorphisms between some of
the $\gR$\hskip -.5mm's in Table 1 given
by \eqref{iso}.

\medskip

Now we restate the classification of
nilpotent ${\rm Ad}(\g^{}_{\R})$-orbits in
all real forms $\gR$ of classical complex
Lie algebras $\g$ (cf.\;\cite{BoCu},
\cite{CM}, \cite{SpSt}, \cite{M},
\cite{W}) in the form adapted to our goal
of classifying compact nilpotent elements
in $\gR$.

If a subalgebra of $\gR$ is isomorphic to
$\ssl_2(\R)$, we call it an
$\ssl_2(\R)$-{\it subalgebra}. We use the
following basic facts, cf.\;\cite{CM},
\cite{KR}, \cite{M}, \cite{SpSt}:

\medskip

(F1) For any nonzero nilpotent element
$x\in\gR$, there is an
$\ssl_2(\R)$-subalgebra of $\gR$
containing~$x$.

(F2) If $\mathfrak a_1$ and $\mathfrak
a_2$ are the $\ssl_2(\R)$-subalgebras of
$\gR$, and a nonzero nilpotent ${\rm
Ad}(\g^{}_{\R})$-orbit intersects both
$\mathfrak a_1$ and $\mathfrak a_2$, then
$\mathfrak a_2=g\cdot \mathfrak a_1$ for
some $g\in {\rm Ad}(\g^{}_{\R})$.
%%Let $\mathfrak a_s$, $s=1, 2$, be
%%$\ssl_2(\R)$-subalgebras of $\gR$ and let
%%$x_s\in \mathfrak a_s$ be a nonzero
%%nilpotent element. If $x_1$ and $x_2$ lie
%%in the same ${\rm Ad}(\g^{}_{\R})$-orbit,
%%then $\mathfrak a_1$ and $\mathfrak a_2$
%%are ${\rm Ad}(\g^{}_{\R})$-conjugate.

(F3) There are exactly two nonzero
nilpotent ${\rm Ad}(\ssl_2(\R))$-orbits in
$\ssl_2(\R)$. Scalar multiplication by
$-1$ maps one of them to the other.

(F4) Any finite dimensional
$\ssl_2(\R)$-$D$-module is completely
reducible. For any integer $d\in \mathbb
N$, there is a unique, up to isomorphism,
$d$-dimensional simple
$\ssl_2(\R)$-$D$-module $S_d$. The module
$S_1$ is trivial. Let $T_m$ be the direct
sum of $m$ copies of $S_1$, and $T_0:=0$.

(F5) If $\mathfrak a$ is an
$\ssl_2(\R)$-subalgebra of $\gR$ and $x\in
\mathfrak a$ is a nonzero nilpotent
element, then
$\z^{}_{\g^{}_{\R}}(\mathfrak a)$ is the
reductive Levi subalgebra of
$\z^{}_{\g^{}_{\R}}\!(x)$.

(F6) If there is a nonzero
$\ssl_2(\R)$-invariant form on $S_d$ of a
given type  (symmetric, skew-symmetric,
Hermitian or skew-Hermitian), it is
nondegenerate and unique up to
proportionality. We fix such a form. Table
2 contains information about its existence
and the notation for the fixed forms.

\

\vskip -8mm

\

\begin{center}{TABLE 2}\\[6pt]
{\fontsize{9pt}{5mm}\selectfont
\begin{tabular}{c|c|cc|c|cc|c}
$D$ & $d$ & \ \hskip 2mm symmetric &
{\hfill \vline} \hskip 1.6mm sgn &
$\begin{matrix}{\rm skew-}\\[-5pt]
{\rm symmetric}\end{matrix}$ &\ \hskip
2mm Hermitian & {\hfill \vline} %%\
\hskip 1.6mm sgn &
$\begin{matrix}{\rm skew-}\\[-5pt]
{\rm Hermitian}\end{matrix}$
\\
 \hline \hline
 &&&&&&&\\[-14pt]
   $\R$
   &
   $\begin{matrix}
\mbox{\rm even}\\[-3pt]
\mbox{\rm odd}
\end{matrix}$
& \ \hskip 2mm
$\begin{matrix}\ \\[-3pt]
\Delta_d^s\end{matrix}$ &
$\begin{matrix}\ \\[-3pt]
1\end{matrix}$\hskip -3mm \ &
$\begin{matrix}\Delta_d^{ss}\\[-3pt] \
\end{matrix}$
& & &
\\[6pt]
\hline
 &&&&&\\[-14pt]
 $\C$
 &
 $\begin{matrix}
\mbox{\rm even}\\[-3pt]
\mbox{\rm odd}
\end{matrix}$
 & \ \hskip 2mm
 $
 \begin{matrix}\
 \\[-3pt] \Delta_d^s\end{matrix}$
 &
& $\begin{matrix} \Delta_d^{ss} \\[-3pt] \
\end{matrix}$
& \ \hskip 2mm $\Delta_d^H$ &
$\begin{matrix}
0\\[-3pt]
1\end{matrix}$ \hskip -3mm \ &
$i\Delta_d^{H}$
\\[6pt]
\hline
 &&&&&&&\\[-14pt]
$\Ha$ & $\begin{matrix}
\mbox{\rm even}\\[-3pt]
\mbox{\rm odd}
\end{matrix}$
 &
 &
 &
 & \ \hskip 2mm
$\begin{matrix}\ \\[-3pt] \Delta_d^H
\end{matrix}$ &
$\begin{matrix}\ \\ 1 \end{matrix}$ \hskip
-3mm \
& $\begin{matrix}\Delta_d^{sH} \\[-3pt]
\ \end{matrix}$
\end{tabular}
}
\end{center}

\vskip 2.5mm

{\it Partitions, $\ssl_2(\R)$-$D$-modules,
and $\ssl_2(\R)$-invariant forms}

\medskip
We call any vector
%%finite sequence
\begin{equation}\label{emm}
\mbox{\boldmath$m$}=(m_1,\ldots,m_p)\in
\mathbb N^p, \ \text{ where } m_p\neq 0,
\end{equation}
%%will be called
a {\it partition} of the number
$$
\textstyle |\mbox{\boldmath$m$}|:=\sum_d
dm_d$$ (this is nontraditional usage of
the term ``partition'' but it is
convenient for our purposes). If
%%$m_1=\ldots=m_p=1$,
$p=1$, then $\mbox{\boldmath$m$}$ is
called a {\it trivial} partition.

The {\it Young diagram} of
\mbox{\boldmath$m$} is a left-justified
array ${\rm Y}(\mbox{\boldmath$m$})$ of
empty boxes with $p$ boxes in each of the
first $m_p$ rows, $p-1$ boxes in each of
the next $m_{p-1}$ rows, and so on. The
partition
\begin{equation*}\label{transp}
\check{\mbox{\boldmath$m$}}=\{{\check
m}_1,\ldots, {\check m}_q\}
\end{equation*}
is called the {\it transpose partition} to
\mbox{\boldmath$m$} if ${\rm
Y}(\check{\mbox{\boldmath$m$}})$ is the
transpose of ${\rm
Y}(\mbox{\boldmath$m$})$, i.e.,\;the rows
of ${\rm Y}(\check{\mbox{\boldmath$m$}})$
are the columns of  ${\rm
Y}(\mbox{\boldmath$m$})$ from left to
right. We have
$|\check{\mbox{\boldmath$m$}}|=
|\mbox{\boldmath$m$}|$.

Let \mbox{\boldmath$m$} be a nontrivial
partition. It is called a {\it symmetric}
(resp.,\;\hskip -.4mm\;{\it
skew-symmetric}) {\it parti\-ti\-on} of
$|\mbox{\boldmath$m$}|$ if $m_d$ in
\eqref{emm} is even for every even
(resp.,\;\hskip -.4mm\;odd) $d$.

Let $\underline{\mbox{\boldmath$m$}}$ be a
sequence obtained from
$\mbox{\boldmath$m$}$ by replacing $m_d$
in \eqref{emm} for every $d$
(resp.,\;\hskip -.4mm\;every odd $d$,
every even $d$) with a pair $(p_d,
q_d)\in\mathbb N^2$ such that
$p_d+q_d=m_d$. Such
$\underline{\mbox{\boldmath$m$}}$ is
called a {\it fine} (resp.,\;\hskip
-.4mm\;{\it fine Hermitian}, {\it fine
skew-Hermitian}) {\it partition} of
$|\mbox{\boldmath$m$}|$ {\it associated}
with $\mbox{\boldmath$m$}$. If
$\underline{\mbox{\boldmath$m$}}$ is fine
or fine Hermitian,
\begin{equation}\label{fsgnn}
\textstyle
 {\rm
sgn}\,\underline{\mbox{\boldmath$m$}}:=
\sum_{d \text{ \rm  odd } }(p_d-q_d)
\end{equation}
is called the {\it signature} of
$\underline{\mbox{\boldmath$m$}}$. If
$\underline{\mbox{\boldmath$m$}}$ is fine
Hermitian (resp.,\;\hskip -.4mm\;fine
skew-Hermitian) and \mbox{\boldmath$m$} is
symmetric (resp.,\;\hskip
-.4mm\;skew-symmetric), then
$\underline{\mbox{\boldmath$m$}}$ is
called a {\it fine symmetric}
(resp.,\;\hskip -.4mm\;{\it fine
skew-symmetric}) {\it partition} of
$|\mbox{\boldmath$m$}|$.

%%\vskip 1.7mm

Any partition \eqref{emm} defines the
$|\mbox{\boldmath$m$}|$-dimensional
$\ssl_2(\R)$-$D$-module
\begin{equation}\label{modulee}
\textstyle V_{\mbox{\boldmath$m$}}:=
\bigoplus_{d\geqslant 1} (T_{m_d}\!\otimes
S_d)
\end{equation}
(tensor product is taken over $D$, with
respect to the canonical $D$-bimodule
structure on $D$-vector spaces). By (F4),
any nonzero finite dimensional
$\ssl_2(\R)$-$D$-module is isomorphic to
$V_{\mbox{\boldmath$m$}}$ for a unique
{\boldmath$m$}. It is trivial if and only
if \mbox{\boldmath$m$} is trivial.

Identifying every $\varphi \in
\gl_D(T_{m_d})$ with the transformation of
$V_{\mbox{\boldmath$m$}}$ acting as
$\varphi\otimes {\rm id}$ on the summand
$T_{m_d}\!\otimes S_d$ and as $0$ on the
other summands in the right-hand side of
\eqref{modulee}, we identify
$\oplus_{d\geqslant 1}\;\gl_D(T_{m_d})$
with the subalgebra of
$\gl_D(V_{\mbox{\boldmath$m$}})$.

Two $\ssl_2(\R)$-invariant forms on
$V_{\mbox{\boldmath$m$}}$ are called
$\ssl_2(\R)$-{\it equivalent} if there is
an $\ssl_2(\R)$-equivariant isometry
$V_{\mbox{\boldmath$m$}}\rightarrow
V_{\mbox{\boldmath$m$}}$ with respect to
them. To describe the
$\ssl_2(\R)$-equivalency classes of
$\ssl_2(\R)$-invariant forms on
$V_{\mbox{\boldmath$m$}}$, we fix, for
every positive integer $r$ and pair $(p,
q)\in \mathbb N^2$ with $p+q=r$, a
nondegenerate form on $T_r$ whose
notation, type and signature (if
applicable) are specified in Table 3.

\

\vskip -8mm

\

\begin{center}{TABLE 3}\\[6pt]
{\fontsize{9pt}{5mm}\selectfont
\begin{tabular}{c|cc|c|cc|c}
$D$ & \hskip 3mm symmetric & {\hfill
\vline} \hskip 2mm sgn
&$\begin{matrix}{\rm skew-}
\\[-5pt]
{\rm symmetric}\end{matrix}$
  &\hskip 3mm
Hermitian &{\hfill \vline} \hskip 2mm sgn
&
$\begin{matrix}{\rm skew-}\\[-5pt]
{\rm Hermitian}\end{matrix}$
\\
 \hline \hline
 &&&&&&\\[-12pt]
   $\R$
   & \hskip 3mm
 $\Theta_{p, q}^s$
 &
  \ $p\!-\!q$\hskip -2.5mm \
  &
$\Theta_r^{ss}$,\ $r$ even & & &
\\[3pt]
\hline
&&&&&&\\[-10pt]
 $\C$
 & \hskip 1.7mm
 $\Theta_r^s$
 &
 &
 $\Theta_r^{ss}$,\ $r$ even &\hskip 4mm
 $\Theta_{p,
q}^H$ & \ $p\!-\!q$ \hskip -2.5mm \ &
\\[3pt]
\hline
&&&&&&\\[-10pt]
$\Ha$
 &
 &
 &
 &\hskip 4mm
 $\Theta_{p, q}^H$
 &
\ $p\!-\!q$ \hskip -2.5mm \ &
$\Theta_r^{sH}$
\end{tabular}
}
\end{center}

\

\vskip -4mm

\

\noindent According to the above
discussion, any nondegenerate form on
$T_r$ is equivalent to a unique form from
Table 3.

If $\Psi$ is a nondegenerate
$\ssl_2(\R)$-invariant form on
$V_{\mbox{\boldmath$m$}}$, then $\Psi|_d$,
its restriction to the summand
$T_{m_d}\!\otimes S_d$ in \eqref{modulee},
is a nondegenerate $\ssl_2(\R)$-invariant
form of the same type (i.e.,\;symmetric,
skew-symmetric, Hermitian or
skew-Hermitian) as $\Psi$, and different
summands in \eqref{modulee} are orthogonal
with respect to $\Psi$. If $\Psi'$ is
another nondegenerate
$\ssl_2(\R)$-invariant form on
$V_{\mbox{\boldmath$m$}}$, then $\Psi$ and
$\Psi'$ are $\ssl_2(\R)$-equivalent if and
only if $\Psi|_d$ and $\Psi'|_d$ are
$\ssl_2(\R)$-equivalent for every $d$. If,
for every $d$, a nondegenerate
$\ssl_2(\R)$-invariant form $\Phi_d$ on
$T_{m_d}\!\otimes S_d$ is fixed, and all
$\Phi_d$'s are of the same type, then
there is a nondegenerate
$\ssl_2(\R)$-invariant form $\Psi$ on
$V_{\mbox{\boldmath$m$}}$ such that
$\Psi|_d=\Phi_d$ for all $d$. This reduces
describing the $\ssl_2(\R)$-equivalency
classes of nondegenerate
$\ssl_2(\R)$-invariant forms on
$V_{\mbox{\boldmath$m$}}$ to describing
that on $T_{m_d}\!\otimes S_d$ for all
$d$. For any positive integers $r$ and
$d$, Table 4 describes all, up to
$\ssl_2(\R)$-equivalency, nondegenerate
$\ssl_2(\R)$-invariant forms on
$T_r\!\otimes S_d$ (in these tables,
 $p+q=r$).

 %%\

 %%\vskip -5mm

 %%\

\begin{center}{TABLE 4}\nopagebreak\\[6pt]
{\fontsize{9pt}{5mm}\selectfont
\begin{tabular}{c|c|cc|c|cc|c}
$D$ & $d$ & \hskip 4mm symmetric & {\hfill
\vline} \hskip 2mm sgn &\hskip 3mm
$\begin{matrix}{\rm skew-}\\[-5pt]
{\rm symmetric}\end{matrix}\hskip 2mm $&
\hskip 3mm Hermitian & {\hfill \vline}
\hskip 2mm sgn&$\hskip 1mm
\begin{matrix}{\rm skew-}\\[-5pt]
{\rm Hermitian}\end{matrix}$
\\
\hline \hline
 &&&&&&&
 \\[-10pt]
   $\R$
   &
   $\begin{matrix}
\mbox{\rm even}\\
\mbox{\rm odd}\end{matrix} $
   &\hskip 5mm
   $\begin{matrix}
 \Theta_r^{ss}\otimes\Delta_{d}^{ss}
 \\
\Theta_{p, q}^s\otimes\Delta_{d}^s
 \end{matrix}$
 &
 $\
 \begin{matrix}0\\
 p\!-\!q\end{matrix}$\hskip -2mm \
&\hskip 1.5mm
 $
 \begin{matrix}
 \Theta_{p, q}^s\otimes\Delta_{d}^{ss}\\
 \Theta_r^{ss}\otimes\Delta_{d}^s
 \end{matrix}$&&&
\\
&&&&&&&
\\[-11pt]
\hline &&&&&&&
\\[-10pt]
 $\C$
 &
$\begin{matrix}
\mbox{\rm even}\\
\mbox{\rm odd}\end{matrix} $
 &\hskip 5mm
 $
\begin{matrix}
 \Theta_r^{ss}\otimes \Delta_{d}^{ss}
 \\
\Theta_r^s\otimes \Delta_{d}^s
 \end{matrix}$
 &
&\hskip 1.5mm $\begin{matrix}\Theta_r^s
\otimes\Delta_{d}^{ss}\\
\Theta_{r}^{ss}\otimes\Delta_{d}^s
\end{matrix}$
& \hskip 3mm $\Theta_{p, q}^H\otimes
\Delta_d^H$ &\
$\begin{matrix}0\\
 p\!-\!q\end{matrix}$\hskip -2mm \
&
\\
&&&&&&&
\\[-10pt]
\hline &&&&&&&
\\[-11pt]
 $\Ha$
 &&&&
 &
 \hskip 3mm
 $\begin{matrix}
 \Theta_{r}^{sH}\otimes\Delta_{d}^{sH}\\
 \Theta_{p, q}^H\otimes
 \Delta_{d}^H
 \end{matrix}$
 &\
 $\begin{matrix}0\\
 p\!-\!q\end{matrix}$\hskip -2mm \
 &
 \hskip 2mm
 $\begin{matrix} \Theta_{p,
q}^H\otimes\Delta_{d}^{sH}\\
\Theta_r^{sH}\otimes\Delta_{d}^H
\end{matrix}$
\end{tabular}
}
\end{center}

\vskip 1.5mm

Returning back to the $n$-dimensional
vector space $V$ over $D$, fix a form
$\Phi$ (not necessarily nondegenerate) on
$V$. Let \mbox{\boldmath$m$} be a
nontrivial partition of $n$ and let
\begin{equation}\label{alpha}
\alpha^{}_{\mbox{\boldmath$m$}}:
\ssl_2(\R)\hookrightarrow
\gl_D(V_{\mbox{\boldmath$m$}})
\end{equation}
be the injection defining the
$\ssl_2(\R)$-$D$-module structure on
$V_{\mbox{\boldmath$m$}}$. Let $\Psi$ be
an $\ssl_2(\R)$-invariant form on
$V_{\mbox{\boldmath$m$}}$. If $\Psi$ and
$\Phi$ are equivalent and $\iota:
V_{\mbox{\boldmath$m$}}\rightarrow V$ is
an isometry with respect to $\Psi$ and
$\Phi$, then the image of the homomorphism
$$\ssl_2(\R)\rightarrow \gl_D(V),\quad
\varphi\mapsto \iota\circ
\alpha^{}_{\mbox{\boldmath$m$}}
(\varphi)\circ \iota^{-1},$$ is an
$\ssl_2(\R)$-subalgebra of $\gl_D^\Phi
(V)$. The above discussion shows that any
$\ssl_2(\R)$-subalgebra of $\gl_D^\Phi
(V)$ is obtained in this way from some
pair $\{\mbox{\boldmath$m$}, \Psi\}$. The
$\ssl_2(\R)$-subalgebras of $\gl_D^\Phi
(V)$ obtained in this way from the pairs
$\{\mbox{\boldmath$m$}, \Psi\}$,
$\{\mbox{\boldmath$m$}', \Psi'\}$ are
${\GL}_D^\Phi (V)$-conjugate if and only
if
$\mbox{\boldmath$m$}=\mbox{\boldmath$m$}'$,
and $\Psi$ and $\Psi'$ are
$\ssl_2(\R)$-equivalent.
 This
yields a bijection between the union of
$\ssl_2(\R)$-equivalency classes of
$\ssl_2(\R)$-invariant forms equivalent to
$\Phi$ on
%%all
$V_{\mbox{\boldmath$m$}}$, where
\mbox{\boldmath$m$} ranges over all
nontrivial partitions of $n$, and the set
of all ${\GL}_D^{\Phi}(V)$-conjugacy
classes of $\ssl_2(\R)$-subalgebras in
$\gl_D^{\Phi}(V)$. {\it Denote this
bijection by} $\star$. If $\Psi$ is either
zero or nondegenerate, considering $\star$
leads to the classification of nilpotent
${\rm Ad}(\gl_D^{\Phi}(V))$-orbits in
$\gl_D^{\Phi}(V)$ described below. In each
case, we give the formulas for the orbit
dimensions (one obtains them using
\eqref{dimKS} and \cite{CM}, \cite{M});
then \eqref{2}, \eqref{dimKS} yield the
dimensions of the corresponding $K$-orbits
in~$\cNp$.

%%\vskip 4mm

\vskip 2.5mm

{\it $\bullet$ Nilpotent orbits in
$\ssl_D(V)$}

\smallskip

Take $\Phi=0$. Then
${\GL}_D^{\Phi}(V)={\GL}_D(V)$ and
$\gl_D^{\Phi}(V)=\gl_D(V)$. Any
$\ssl_2(\R)$-subalgebra of $\gl_D(V)$ is
contained in $\ssl_D(V)$. For any
nontrivial partition {\boldmath$m$} of
$n$, consider the ${\rm
GL}_D^{\Phi}(V)$-conjugacy class of
$\ssl_2(\R)$-subalgebras in $\ssl_D(V)$
corresponding under $\star$ to the
$\ssl_2(\R)$-equivalency class of zero
form  on $V_{\mbox{\boldmath$m$}}$. Then
by (F1)--(F3), \cite{KR}, the subset of
$\gl_D^{\Phi}(V)$ that is the union of all
$\ssl_2(\R)$-subalgebras from this class
contains a unique nonzero nilpotent
${\GL}_D(V)$-orbit~${\cal
O}_{\mbox{\boldmath$m$}}$.

\smallskip

${\rm Ad}(\ssl_D(V))$-orbits in ${\cal
N}(\ssl_D(V))$:
\begin{quotation}\noindent
${\cal O}_{\mbox{\boldmath$m$}}$ is an
${\rm Ad}(\ssl_D(V))$-orbit,  except that
${\cal O}_{\mbox{\boldmath$m$}}$ is the
union of two ${\rm Ad}(\ssl_D(V))$-orbits
${\cal O}_{\mbox{\boldmath$m$}}^1$ and
${\cal O}_{\mbox{\boldmath$m$}}^2$ (that
are the connected components of ${\cal
O}_{\mbox{\boldmath$m$}}$) if $D=\R$ and
$m_d=0$ for every odd $d$. %%
Such ${\rm Ad}(\ssl_D(V))$-orbits, taken
over all partitions \mbox{\boldmath$m$} of
$n$, are pairwise different and exhaust
all nonzero
 ${\rm Ad}(\ssl_D(V))$-orbits in
${\cal N}(\ssl_D(V))$.
\end{quotation}

Orbit dimensions:
%% (see
%%\eqref{emm}, \eqref{transp}):
\begin{gather*}\textstyle
n^2-\sum_d d^2{\check m}_d =\begin{cases}
\dim^{}_{\R}{\cal O}_{\mbox{\boldmath$m$}}
&\text{\rm for }
D=\R,\\
 \dim^{}_{\C}{\cal
O}_{\mbox{\boldmath$m$}}
&\text{\rm for } D=\C,\\
%%\frac{1}{4}
\dim^{}_{\R}{\cal
O}_{\mbox{\boldmath$m$}}/4 &\text{\rm for
} D=\Ha.
\end{cases}
%%\\
%%\textstyle 4n^2-4\sum_d d^2{\check m}_d
%%=\dim^{}_{\R}{\cal
%%O}_{\mbox{\boldmath$m$}}\ \text{\rm for }
%%D=\Ha.
\end{gather*}

%%\vskip 2mm

{\it  $\bullet$ Nilpotent orbits in
$\gl_{D}^\Phi (V)$ for $D=\R$ and
nondegenerate symmetric $\Phi$}

\smallskip

Let $\underline{\mbox{\boldmath$m$}}$ be a
fine symmetric partition of $n$ associated
with \eqref{emm}.
%%${\mbox{\boldmath$m$}}$
%%the associated
%%partition.
Let
$\Psi_{\underline{\mbox{\boldmath$m$}}}$
be the $\ssl_2(\R)$-invariant form on
$V_{\mbox{\boldmath$m$}}$ such that for
all $d$,
\begin{equation}\label{symmR}
\Psi_{\underline{\mbox{\boldmath$m$}}}|_d=
\begin{cases}
\Theta_{p_d, q_d}^s\otimes \Delta_d^s
&\text{ \rm if } d \text{ \rm is odd, }\\
 \Theta_{m_d}^{ss}\otimes \Delta_d^{ss}
&\text{ \rm if } d \text{ \rm is even. }
\end{cases}
\end{equation}
 Then
$\Psi_{\underline{\mbox{\boldmath$m$}}}$
is equivalent to $\Phi$ if and only if
\begin{equation}\label{sgn}{\rm
sgn}\;{\underline{\mbox{\boldmath$m$}}}
={\rm sgn}\,\Phi.\end{equation} If
\eqref{sgn} holds, consider  the ${\rm
GL}_D^{\Phi}(V)$-conjugacy class of
$\ssl_2(\R)$-subalgebras in
$\gl_D^{\Phi}(V)$  corresponding under
%%the bijection
$\star$ to the $\ssl_2(\R)$-equivalency
class of
$\Psi_{\underline{\mbox{\boldmath$m$}}}$.
Then by (F1)--(F3), \cite{KR}, the union
of all $\ssl_2(\R)$-subalgebras from this
class contains a unique nilpotent ${\rm
GL}_D^{\Phi}(V)$-orbit~${\cal
O}_{\underline{\mbox{\boldmath$m$}}}$.

\smallskip

${\rm Ad}(\gl_D^{\Phi}(V))$-orbits in
${\cal N}(\gl_D^{\Phi}(V))$:
\begin{quotation}

\noindent If \eqref{sgn} holds, ${\cal
O}_{\underline{\mbox{\boldmath$m$}}}$ is
an ${\rm Ad}(\gl_D^{\Phi}(V))$-orbit,
except the following cases. If $m_d=0$ for
all odd $d$, then ${\cal
O}_{\underline{\mbox{\boldmath$m$}}}$ is
the union of four ${\rm
Ad}(\gl_D^{\Phi}(V))$-orbits ${\cal
O}_{\underline{\mbox{\boldmath$m$}}}^1$,
${\cal
O}_{\underline{\mbox{\boldmath$m$}}}^2$,
${\cal
O}_{\underline{\mbox{\boldmath$m$}}}^3$,
${\cal
O}_{\underline{\mbox{\boldmath$m$}}}^4$
that are the connected components of
${\cal
O}_{\underline{\mbox{\boldmath$m$}}}$. If
there is an odd $d$ such that $m_d\neq 0$
and
$$
\mbox{\text{\rm either }}\hskip 1mm
\begin{cases}
p_d=0 &\mbox{\text{\rm for all }}
d\equiv 1\,{\rm mod}\,4,\\
q_d=0 &\mbox{\text{\rm for all }} d\equiv
3\,{\rm mod}\,4
\end{cases}\hskip 3mm
\mbox{\text{\rm or }}\hskip 1mm
\begin{cases} p_d=0
&\mbox{\text{\rm for all }}
d\equiv 3\,{\rm mod}\,4,\\
 q_d=0 &\mbox{\text{\rm for all }}
d\equiv 1\,{\rm mod}\,4,
\end{cases}
$$
%%and $p_dq_d=0$ for all such $d$,
then ${\cal
O}_{\underline{\mbox{\boldmath$m$}}}$ is
the union of two ${\rm
Ad}(\gl_D^{\Phi}(V))$-orbits ${\cal
O}_{\underline{\mbox{\boldmath$m$}}}^1$,
${\cal
O}_{\underline{\mbox{\boldmath$m$}}}^2$
that are the connected components of
${\cal
O}_{\underline{\mbox{\boldmath$m$}}}$.
Such ${\rm Ad}(\gl_D^{\Phi}(V))$-orbits,
taken over all fine symmetric partitions
$\underline{\mbox{\boldmath$m$}}$ of $n$
satisfying \eqref{sgn}, are pairwise
different  and  exhaust all nonzero ${\rm
Ad}(\gl_D^{\Phi}(V))$-orbits in ${\cal
N}(\gl_D^{\Phi}(V))$.
\end{quotation}

Orbit dimensions:
%%(see
%%\eqref{emm},\eqref{transp}):
\begin{equation*}
\textstyle \dim^{}_{\R}{\cal
O}_{\underline{\mbox{\boldmath$m$}}}=
%%\frac{1}{2}
\bigl(n^2-n-\sum_d d^2 {\check m}_d
+\sum_{d \text{ \rm odd }} m_d\bigr)/2.
\end{equation*}

%%\vskip 2mm

{\it $\bullet$ Nilpotent orbits in
$\gl_{D}^\Phi (V)$ for $D=\C$ and
nondegenerate symmetric $\Phi$}

\smallskip

Let $\mbox{\boldmath$m$}$ be a symmetric
partition of $n$. Let
$\Psi_{{\mbox{\boldmath$m$}}}$ be the
$\ssl_2(\R)$-invariant form on
$V_{\mbox{\boldmath$m$}}$ such that for
all $d$,
$$
\Psi_{\mbox{\boldmath$m$}}|_d=
\begin{cases}
\Theta_{m_d}^{s}\otimes \Delta_d^s
&\text{ \rm if } d \text{ \rm is odd, }\\
 \Theta_{m_d}^{ss}\otimes \Delta_d^{ss}
 &\text{ \rm if } d \text{ \rm is even.
}
\end{cases}
$$
Then $\Psi_{{\mbox{\boldmath$m$}}}$ is
equivalent to $\Phi$. Consider  the ${\rm
GL}_D^{\Phi}(V)$-conjugacy class of
$\ssl_2(\R)$-subalgebras in
$\gl_D^{\Phi}(V)$ corresponding under
$\star$ to the $\ssl_2(\R)$-equivalency
class of $\Psi_{{\mbox{\boldmath$m$}}}$.
Then by (F1)--(F3), \cite{KR}, the union
of all $\ssl_2(\R)$-subalgebras from this
class contains a unique nilpotent ${\rm
GL}_D^{\Phi}(V)$-orbit~${\cal
O}_{{\mbox{\boldmath$m$}}}$.

\smallskip

${\rm Ad}(\gl_D^{\Phi}(V))$-orbits in
${\cal N}(\gl_D^{\Phi}(V))$:
\begin{quotation}\noindent
${\cal O}_{\mbox{\boldmath$m$}}$ is an
${\rm Ad}(\gl_D^{\Phi}(V))$-orbit, except
that if $m_d=0$ for all odd $d$, then
${\cal O}_{{\mbox{\boldmath$m$}}}$ is the
union of two ${\rm
Ad}(\gl_D^{\Phi}(V))$-orbits  ${\cal
O}_{{\mbox{\boldmath$m$}}}^1$ and ${\cal
O}_{{\mbox{\boldmath$m$}}}^2$ that are the
connected components of ${\cal
O}_{\mbox{\boldmath$m$}}$. Such ${\rm
Ad}(\gl_D^{\Phi}(V))$-orbits, taken over
all symmetric par\-titions
 \mbox{\boldmath$m$} of $n$,
are pairwise different and exhaust all
nonzero ${\rm Ad}(\gl_D^{\Phi}(V))$-orbits
in ${\cal N}(\gl_D^{\Phi}(V))$.
\end{quotation}

Orbit dimensions:
%%(see
%%\eqref{emm},\eqref{transp}):
\begin{equation*}
\textstyle \dim^{}_{\C}{\cal
O}_{\mbox{\boldmath$m$}}=
%%\frac{1}{2}
\bigl(n^2-n-\sum_d d^2 {\check m}_d
+\sum_{d \text{ \rm odd }} m_d\bigr)/2.
\end{equation*}

%%\vskip 2mm

{\it $\bullet$ Nilpotent orbits in
$\gl_{D}^\Phi (V)$ for $D=\R$ and
nondegenerate skew-symmetric $\Phi$}

\smallskip

Let $\underline{\mbox{\boldmath$m$}}$ be a
fine skew-symmetric partition of $n$
associated with \eqref{emm}.
 Let
$\Psi_{\underline{\mbox{\boldmath$m$}}}$
be the $\ssl_2(\R)$-invariant form on
$V_{\mbox{\boldmath$m$}}$ such that for
all $d$,
\begin{equation}\label{skewR}
\Psi_{\underline{\mbox{\boldmath$m$}}}|_d=
\begin{cases}
\Theta_{p_d, q_d}^s\otimes \Delta_d^{ss}
&\text{ \rm if } d \text{ \rm is even, }\\
 \Theta_{m_d}^{ss}\otimes \Delta_d^s
&\text{ \rm if } d \text{ \rm is odd. }
\end{cases}
\end{equation}
Then
$\Psi_{\underline{\mbox{\boldmath$m$}}}$
is equivalent to $\Phi$. Consider  the
${\rm GL}_D^{\Phi}(V)$-conjugacy class of
$\ssl_2(\R)$-subalgebras in
$\gl_D^{\Phi}(V)$ corresponding under
$\star$ to the $\ssl_2(\R)$-equivalency
class of
$\Psi_{\underline{\mbox{\boldmath$m$}}}$.
Then by (F1)--(F3), \cite{KR}, the union
of all $\ssl_2(\R)$-subalgebras from this
class contains a unique nilpotent ${\rm
GL}_D^{\Phi}(V)$-orbit ${\cal
O}_{\underline{\mbox{\boldmath$m$}}}$.

\smallskip

${\rm Ad}(\gl_D^{\Phi}(V))$-orbits in
${\cal N}(\gl_D^{\Phi}(V))$:
\begin{quotation}\noindent
${\cal
O}_{\underline{\mbox{\boldmath$m$}}}$ is
an ${\rm Ad}(\gl_D^{\Phi}(V))$-orbit. Such
${\rm Ad}(\gl_D^{\Phi}(V))$-orbits, taken
over all fine skew-symmet\-ric partitions
\underline{\mbox{\boldmath$m$}} of $n$,
are pairwise different  and exhaust all
 nonzero  ${\rm
Ad}(\gl_D^{\Phi}(V))$-orbits in ${\cal
N}(\gl_D^{\Phi}(V))$.
\end{quotation}

Orbit dimensions:
%% (see
%%\eqref{emm}, \eqref{transp}):
\begin{equation*}
\textstyle \dim^{}_{\R} {\cal
O}_{\underline{\mbox{\boldmath$m$}}}=
%%\frac{1}{2}
\bigl(n^2+n-\sum_d d^2 {\check m}_d
-\sum_{d \text{ \rm odd }} m_d\bigr)/2.
\end{equation*}

%%\vskip 4mm

{\it $\bullet$ Nilpotent orbits in
$\gl_{D}^\Phi (V)$ for $D=\C$ and
nondegenerate skew-symmetric $\Phi$}

\smallskip

Let \mbox{\boldmath$m$} be a
skew-symmetric partition of $n$.
 Let
$\Psi_{\mbox{\boldmath$m$}}$ be the
$\ssl_2(\R)$-invariant form on
$V_{\mbox{\boldmath$m$}}$ such that for
all $d$,
$$
\Psi_{\underline{\mbox{\boldmath$m$}}}|_d=
\begin{cases}
\Theta_{m_d}^{s}\otimes \Delta_d^{ss}
&\text{ \rm if } d \text{ \rm is even, }\\
\Theta_{m_d}^{ss}\otimes \Delta_d^s
&\text{ \rm if } d \text{ \rm is odd. }
\end{cases}
$$
Then
$\Psi_{\underline{\mbox{\boldmath$m$}}}$
is equivalent to $\Phi$. Consider  the
${\rm GL}_D^{\Phi}(V)$-conjugacy class of
$\ssl_2(\R)$-subalgebras in
$\gl_D^{\Phi}(V)$ corresponding under
 $\star$ to the
$\ssl_2(\R)$-equivalency class of
$\Psi_{\mbox{\boldmath$m$}}$. Then by
(F1)--(F3), \cite{KR}, the union of all
$\ssl_2(\R)$-subalgebras from this class
contains a unique nilpotent ${\rm
GL}_D^{\Phi}(V)$-orbit ${\cal
O}_{\mbox{\boldmath$m$}}$.

\smallskip

 ${\rm
Ad}(\gl_D^{\Phi}(V))$-orbits in ${\cal
N}(\gl_D^{\Phi}(V))$:
\begin{quotation}\noindent
${\cal O}_{\mbox{\boldmath$m$}}$ is an
${\rm Ad}(\gl_D^{\Phi}(V))$-orbit. Such
${\rm Ad}(\gl_D^{\Phi}(V))$-orbits, taken
over all skew-sym\-metric partitions
\mbox{\boldmath$m$} of $n$, are pairwise
different  and exhaust all nonzero ${\rm
Ad}(\gl_D^{\Phi}(V))$-orbits in ${\cal
N}(\gl_D^{\Phi}(V))$.
\end{quotation}

Orbit dimensions:
%% (see
%%\eqref{emm}, \eqref{transp}):
\begin{equation*}
\textstyle \dim^{}_{\C} {\cal
O}_{\mbox{\boldmath$m$}}=
%%\frac{1}{2}
\bigl(n^2+n-\sum_d d^2 {\check m}_d
-\sum_{d \text{ \rm odd }} m_d\bigr)/2.
\end{equation*}

%%\vskip 4mm

{\it $\bullet$ Nilpotent orbits in
$\gl_{D}^\Phi (V)$ for $D=\C$ and
nondegenerate Hermitian $\Phi$}

\smallskip

Let $\underline{\mbox{\boldmath$m$}}$ be a
fine partition of $n$ associated with
\eqref{emm}. Let
$\Psi_{\underline{\mbox{\boldmath$m$}}}$
be the $\ssl_2(\R)$-invariant form on
$V_{\mbox{\boldmath$m$}}$ such that for
all $d$,
\begin{equation}\label{HermC}
\Psi_{\underline{\mbox{\boldmath$m$}}}|_d=
\Theta_{p_d, q_d}^H\otimes \Delta_d^H.
\end{equation}
Then
$\Psi_{\underline{\mbox{\boldmath$m$}}}$
is equivalent to $\Phi$ if and only if
\eqref{sgn} holds. In the last case,
consider the ${\rm
GL}_D^{\Phi}(V)$-conjugacy class of
$\ssl_2(\R)$-subalgebras in
$\gl_D^{\Phi}(V)$ corresponding under
$\star$ to the $\ssl_2(\R)$-equivalency
class of
$\Psi_{\underline{\mbox{\boldmath$m$}}}$.
Then by (F1)--(F3), \cite{KR}, the union
of all $\ssl_2(\R)$-subalgebras from this
class contains a unique nilpotent ${\rm
GL}_D^{\Phi}(V)$-orbit ${\cal
O}_{\underline{\mbox{\boldmath$m$}}}$.

\smallskip

 ${\rm
Ad}(\gl_D^{\Phi}(V))$-orbits in ${\cal
N}(\gl_D^{\Phi}(V))$:

\begin{quotation}\noindent
 If \eqref{sgn} holds, ${\cal
O}_{\underline{\mbox{\boldmath$m$}}}$ is
an ${\rm Ad}(\gl_D^{\Phi}(V))$-orbit. Such
${\rm Ad}(\gl_D^{\Phi}(V))$-orbits, taken
over all fine partitions
$\underline{\mbox{\boldmath$m$}}$ of $n$
sa\-tisfying \eqref{sgn}, are pairwise
different and exhaust all nonzero ${\rm
Ad}(\gl_D^{\Phi}(V))$-orbits in ${\cal
N}(\gl_D^{\Phi}(V))$.
\end{quotation}

Orbit dimensions:
%% (see
%%\eqref{emm}, \eqref{transp}):
\begin{equation*}
\dim^{}_{\R}{\cal
O}_{\underline{\mbox{\boldmath$m$}}}=\textstyle
n^2-\sum_d d^2{\check m}_d.
\end{equation*}

%%\vskip 2mm

{\it $\bullet$ Nilpotent orbits in
$\gl_{D}^\Phi (V)$ for $D=\Ha$ and
nondegenerate Hermitian $\Phi$}

\smallskip

Let $\underline{\mbox{\boldmath$m$}}$ be a
fine Hermitian partition of $n$ associated
with \eqref{emm}.
 Let
$\Psi_{\underline{\mbox{\boldmath$m$}}}$
be the $\ssl_2(\R)$-invariant form on
$V_{\mbox{\boldmath$m$}}$ such that for
all $d$,
\begin{equation}\label{HermH}
\Psi_{\underline{\mbox{\boldmath$m$}}}|_d=
\begin{cases}
\Theta_{m_d}^{sH}\otimes \Delta_d^{sH}
&\text{ \rm if } d \text{ \rm is even, }\\
 \Theta_{p_d, q_d}^H\otimes \Delta_d^H
&\text{ \rm if } d \text{ \rm is odd. }
\end{cases}
\end{equation}
Then
$\Psi_{\underline{\mbox{\boldmath$m$}}}$
is equivalent to $\Phi$ if and only if
\eqref{sgn} holds. In the last case,
consider  the ${\rm
GL}_D^{\Phi}(V)$-conjugacy class of
$\ssl_2(\R)$-subalgebras in
$\gl_D^{\Phi}(V)$ corresponding under
$\star$ to the $\ssl_2(\R)$-equivalency
class of
$\Psi_{\underline{\mbox{\boldmath$m$}}}$.
Then by (F1)--(F3), \cite{KR}, the union
of all $\ssl_2(\R)$-subalgebras from this
class contains a unique nilpotent ${\rm
GL}_D^{\Phi}(V)$-orbit ${\cal
O}_{\underline{\mbox{\boldmath$m$}}}$.

\smallskip

 ${\rm
Ad}(\gl_D^{\Phi}(V))$-orbits in ${\cal
N}(\gl_D^{\Phi}(V))$:
\begin{quotation}
\noindent If \eqref{sgn} holds, ${\cal
O}_{\underline{\mbox{\boldmath$m$}}}$ is
an ${\rm Ad}(\gl_D^{\Phi}(V))$-orbit. Such
${\rm Ad}(\gl_D^{\Phi}(V))$-orbits, taken
over all fine Hermitian partitions
\underline{\mbox{\boldmath$m$}} of $n$
sa\-tisfying \eqref{sgn}, are pairwise
different  and exhaust all nonzero ${\rm
Ad}(\gl_D^{\Phi}(V))$-orbits in ${\cal
N}(\gl_D^{\Phi}(V))$.
\end{quotation}

Orbit dimensions:
%% (see
%%\eqref{emm}, \eqref{transp}):
\begin{equation*}\textstyle
\dim^{}_{\R}{\cal
O}_{\underline{\mbox{\boldmath$m$}}}=
2n^2+n-2\sum_d d^2{\check m}_d-\sum_{d
\text{ \rm odd }} m_d.
\end{equation*}

%%\vskip 2mm

{\it $\bullet$ Nilpotent orbits in
$\gl_{D}^\Phi (V)$ for $D=\Ha$ and
skew-Hermitian $\Phi$}

\smallskip

Let $\underline{\mbox{\boldmath$m$}}$ be a
fine skew-Hermitian partition of $n$
associated with \eqref{emm}.
 Let
$\Psi_{\underline{\mbox{\boldmath$m$}}}$
be the $\ssl_2(\R)$-invariant form on
$V_{\mbox{\boldmath$m$}}$ such that for
all $d$,
\begin{equation}\label{sHermH}
\Psi_{\underline{\mbox{\boldmath$m$}}}|_d=
\begin{cases}
\Theta_{p_d, q_d}^{H}\otimes \Delta_d^{sH}
&\text{ \rm if } d \text{ \rm is even, }\\
 \Theta_{m_d}^{sH}\otimes \Delta_d^H
&\text{ \rm if } d \text{ \rm is odd. }
\end{cases}
\end{equation}
Then
$\Psi_{\underline{\mbox{\boldmath$m$}}}$
is equivalent to $\Phi$. Consider  the
${\rm GL}_D^{\Phi}(V)$-conjugacy class of
$\ssl_2(\R)$-subalgebras in
$\gl_D^{\Phi}(V)$ corresponding under
$\star$ to the $\ssl_2(\R)$-equivalency
class of
$\Psi_{\underline{\mbox{\boldmath$m$}}}$.
Then by (F1)--(F3), \cite{KR}, the union
of all $\ssl_2(\R)$-subalgebras from this
class contains a unique nilpotent ${\rm
GL}_D^{\Phi}(V)$-orbit ${\cal
O}_{\underline{\mbox{\boldmath$m$}}}$.

\smallskip

 ${\rm
Ad}(\gl_D^{\Phi}(V))$-orbits in ${\cal
N}(\gl_D^{\Phi}(V))$:
\begin{quotation}\noindent
 ${\cal
O}_{\underline{\mbox{\boldmath$m$}}}$ is
an ${\rm Ad}(\gl_D^{\Phi}(V))$-orbit. Such
${\rm Ad}(\gl_D^{\Phi}(V))$-orbits, taken
over all fine skew-Hermitian par\-titions
\underline{\mbox{\boldmath$m$}} of $n$,
are pairwise different  and exhaust all
${\rm Ad}(\gl_D^{\Phi}(V))$-orbits in
${\cal N}(\gl_D^{\Phi}(V)\!)$.
\end{quotation}

Orbit dimensions:
%% (see
%%\eqref{emm},\eqref{transp}):
\begin{equation*}\textstyle
\dim^{}_{\R}{\cal
O}_{\underline{\mbox{\boldmath$m$}}}=
2n^2-n-2\sum_d d^2{\check m}_d+\sum_{d
\text{ \rm odd }} m_d.
\end{equation*}

\smallskip

Now we are ready to classify compact
nilpotent elements in the real forms of
complex classical Lie algebras. Recall
that $n=\dim_D V$.

%%\vskip 4mm

\vskip 2mm

{\it $\bullet$ Compact nilpotent elements
in $\ssl_n(\R)$ and $\ssl_n(\Ha)$}

\begin{theorem}\label{sl-compa}
Let ${\cal O}$ be the ${\rm
Ad}(\ssl_D(V))$-orbit of a nonzero element
$x\in{\cal N}(\ssl_D(V))$ where $D$ is
$\R$ or $\Ha$. Let $\mathfrak a$ be an
$\ssl_2(\R)$-subalgebra of $\ssl_D(V)$
containing $x$. The following properties
are equivalent: \blist

\item[\rm (i)] $x$ is compact.

\item[\rm(ii)] $\mathfrak a$-module $V$ is
simple.

\item[\rm(iii)] If $D=\R$, then ${\cal
O}={\cal O}_{(0,\ldots ,0, 1)}$ for odd
$n$, and ${\cal O}_{(0,\ldots ,0, 1)}^1$
or ${\cal O}_{(0,\ldots , 0, 1)}^2$ for
even $n$. If $D=\Ha$, then ${\cal O}={\cal
O}_{(0,\ldots ,0, 1)}$.

\end{list}
\end{theorem}

\noindent {\it Proof.} According to the
previous discussion, we may, and will,
assume that $V=V_{\mbox{\boldmath$m$}}$
and $\mathfrak
a=\alpha_{\mbox{\boldmath$m$}}(\ssl_2(\R))$
(see \eqref{alpha}) for some nontrivial
partition \mbox{\boldmath$m$} of $n$. The
Double Centralizer Theorem  implies that
\begin{equation}\label{zentr}\textstyle
\z_{\gl_D(V)}(\mathfrak
a)=\bigoplus_{d\geqslant
1}\;\gl_D(T_{m_d}).
\end{equation}
Since $\gl_D(T_{m_d})$ (resp.,\;\hskip
-.4mm\;$\ssl_D(T_{m_d})$) is compact if
and only if $m_d=0$ (resp.,\;\hskip
-.4mm\;$0$ or $1$), the claim follows from
\eqref{zentr} and (F5). \quad $\square$

%%\pagebreak

%%\vskip 4mm

%%\nopagebreak

\vskip 2mm

{\it $\bullet$ Compact nilpotent elements
in $\so_{n-q, q}$}

\begin{theorem}\label{symR-compa}
Let ${\cal O}$ be the ${\rm
Ad}(\gl_{\R}^{\Phi}(V))$-orbit of a
nonzero element $x\in{\cal
N}(\gl_{\R}^\Phi (V))$ for a
nondege\-ne\-rate symmetric form $\Phi$.
The following properties are equivalent:
\blist

\item[\rm (i)] $x$ is compact.

\item[\rm(ii)] ${\cal O}= {\cal
O}_{\underline{\mbox{\boldmath$m$}}}^1$ or
${\cal
O}_{\underline{\mbox{\boldmath$m$}}}^2$,
where $\underline{\mbox{\boldmath$m$}}$ is
a fine symmetric partition
 of $n$
such that \eqref{sgn} holds, $p_dq_d=0$
for all odd $d$, and $m_d=0$ for all even
$d$.
\end{list}
\end{theorem}

\noindent {\it Proof.} We may, and will,
assume that $V=V_{\mbox{\boldmath$m$}}$,
$\Phi$ is
$\Psi_{\underline{\mbox{\boldmath$m$}}}$
defined by \eqref{symmR}, and
$x\in\mathfrak
a:=\alpha_{\mbox{\boldmath$m$}}
(\ssl_2(\R))$ for some fine symmetric
partition
$\underline{\mbox{\boldmath$m$}}$ of $n$
such that \eqref{sgn} holds. Then
\eqref{zentr} holds, whence by
\eqref{symmR}
\begin{equation}\label{zentrF1}\textstyle
\z_{\gl_\R^{\Phi}(V)}(\mathfrak
a)=\bigoplus_{d\geqslant
1}\;\gl_\R^{\Theta_d}(T_{m_d})\ \mbox{ \rm
where } \ \Theta_d=\begin{cases}
\Theta_{p_d, q_d}^s &\mbox{ \rm if $d$ is
odd},\\  \Theta_{m_d}^{ss}&\mbox{ \rm if
$d$ is even.}
\end{cases}
\end{equation}
Since $\gl_\R^{\Theta_d}(T_{m_d})$ for
$\Theta_d$ given by \eqref{zentrF1}
%%$\Theta_{p_d, q_d}^s$
%%(resp.,\;\hskip
%%-.4mm\;$\Theta_{m_d}^{ss}$)
is compact if and only if $p_dq_d=0$ for
odd $d$ and $m_d=0$ for even
$d$
%%(resp.,\;\hskip -.4mm\;$m_d=0$)
(see
Table 1 and \eqref{iso}), the claim
follows from \eqref{zentrF1} and
(F5).\quad $\square$

%%\vskip 4mm
\vskip 2mm

{\it $\bullet$ Compact nilpotent elements
in $\spl_{n}(\R)$}

\begin{theorem}\label{skew-symR-compa}
Let ${\cal O}$ be the ${\rm
Ad}(\gl_{\R}^{\Phi}(V))$-orbit of a
nonzero element $x\in{\cal
N}(\gl_{\R}^\Phi (V))$ for a
nondege\-ne\-rate skew-symmetric form
$\Phi$. The following properties are
equivalent: \blist

\item[\rm (i)] $x$ is compact.

\item[\rm(ii)] ${\cal O}={\cal
O}_{\underline{\mbox{\boldmath$m$}}}$
where $\underline{\mbox{\boldmath$m$}}$ is
a fine skew-symmetric partition
 of $n$
such that $p_dq_d=0$ for all even $d$, and
$m_d=0$ for all odd $d$.
\end{list}
\end{theorem}
\noindent {\it Proof.} The arguments are
similar to that in the proof of Theorem
\ref{symR-compa} with \eqref{symmR}
replaced with \eqref{skewR}. \quad
$\square$

%%\vskip 4mm
\vskip 2mm

{\it $\bullet$ Compact nilpotent elements
in $\su_{n-q, q}$}

\begin{theorem}\label{HermC-compa}
Let ${\cal O}$ be the ${\rm
Ad}(\gl_{\C}^{\Phi}(V))$-orbit of a
nonzero element $x\in{\cal
N}(\gl_{\C}^\Phi (V))$ for a
nondege\-ne\-rate Hermitian form $\Phi$.
The following properties are equivalent:
\blist

\item[\rm (i)] $x$ is compact.

\item[\rm(ii)] ${\cal O}={\cal
O}_{\underline{\mbox{\boldmath$m$}}}$
where $\underline{\mbox{\boldmath$m$}}$ is
a fine partition
 of $n$
such that \eqref{sgn} holds and $p_dq_d=0$
for all $d$.
\end{list}
\end{theorem}
\noindent {\it Proof.} By Table 1, the
algebra
%%the
%%form $\Theta_{p_d, q_d}^H$
$\gl_{\C}^{\Psi}(W)$ for
$\Psi=\Theta_{p_d, q_d}^H$ and $W=T_{m_d}$
is compact if and only if $p_dq_d=0$. Now
one completes the proof using the
arguments similar to that in the proof of
Theorem \ref{symR-compa} with
\eqref{symmR} replaced with \eqref{HermC}.
\quad $\square$

%%\vskip 4mm
\vskip 2mm

{\it $\bullet$ Compact nilpotent elements
in $\spl_{n/2-q, q}$}

\begin{theorem}\label{HermH-compa}
Let ${\cal O}$ be the ${\rm
Ad}(\gl_{\Ha}^{\Phi}(V))$-orbit of a
nonzero element $x\in{\cal
N}(\gl_{\Ha}^\Phi (V))$ for a
nondege\-ne\-rate Hermitian form $\Phi$.
The following properties are equivalent:
\blist

\item[\rm (i)] $x$ is compact.

\item[\rm(ii)] ${\cal O}={\cal
O}_{\underline{\mbox{\boldmath$m$}}}$
where $\underline{\mbox{\boldmath$m$}}$ is
a fine Hermitian partition
 of $n$
such that \eqref{sgn} holds, $p_dq_d=0$
for all odd $d$, and $m_d=0$ or $1$ for
all even $d$.
\end{list}
\end{theorem}
\noindent {\it Proof.} It follows from
Table 1 and \eqref{iso} that the algebra
$\gl_{\Ha}^{\Psi}(W)$ for
$\Psi=\Theta_{p_d, q_d}^H$ (resp.,\;\hskip
-.4mm\;$\Theta_{m_d}^{sH}$) and
$W=T_{m_d}$
%%$\Theta_{p_d, q_d}^H$ (resp.,\;\hskip
%%-.4mm\;$\Theta_{m_d}^{sH}$)
is compact if
and only if $p_dq_d=0$ (resp.,\;\hskip
-.4mm\;$m_d=0$ or $1$). Now one completes
the proof using the arguments similar to
that in the proof of Theorem
\ref{symR-compa} with \eqref{symmR}
replaced with \eqref{HermH}. \quad
$\square$

%%\vskip 4mm
%%\vskip 2mm
\pagebreak

{\it $\bullet$ Compact nilpotent elements
in $\un_{n}^*(\Ha)$}

%%\nopagebreak

\begin{theorem}\label{skew-HermH-compa}
Let ${\cal O}$ be the ${\rm
Ad}(\gl_{\Ha}^{\Phi}(V))$-orbit of a
nonzero element $x\in{\cal
N}(\gl_{\Ha}^\Phi (V))$ for a
nondege\-ne\-rate skew-Hermitian form
$\Phi$. The following properties are
equivalent: \blist

\item[\rm (i)] $x$ is compact.

\item[\rm(ii)] ${\cal O}={\cal
O}_{\underline{\mbox{\boldmath$m$}}}$
where $\underline{\mbox{\boldmath$m$}}$ is
a fine skew-Hermitian partition
 of $n$
such that $p_dq_d=0$ for all even $d$, and
$m_d=0$ or $1$ for all odd $d$.
\end{list}
\end{theorem}
\noindent {\it Proof.}  The arguments are
similar to that in the proof of Theorem
\ref{HermH-compa} with \eqref{HermH}
replaced with \eqref{sHermH}. \quad
$\square$

 \section*{5. Classification for
 exceptional
 simple \boldmath$\g$
 \label{exept}}

In this section, we assume that $\g$ is an
exceptional simple complex Lie algebra and
$\gR$ is its noncompact real form.
According to E.\;Cartan's classification
(cf.,\,e.g.,\,\cite{OV}), up to
isomorphism there are, in total, twelve
possibilities listed in Table~5. Its
second column gives both types of
E.\;Cartan's notation for $\gR$; the
notation ${\mbox{\tt
\fontsize{13pt}{0mm}\selectfont
X}}_{s(t)}$ means that ${\mbox{\tt
\fontsize{13pt}{0mm}\selectfont X}}_s$ is
the type of $\g$ and $t$ is the signature
of the Killing form of $\gR$, i.e.,
$t=\dim^{}_{\C} \pe -\dim^{}_{\C} \ka$.

\vskip 3mm

\begin{center}{TABLE 5}\\[9pt]
{\fontsize{9pt}{5mm}\selectfont
\begin{tabular}{c|c|c|c|c}
\ $\hskip .2mm \g\hskip 1.5mm $ & $\gR$ &
type of $\kR$ &$\dim^{}_{\C} \ka$
&$\dim^{}_{\C} \p$ \\
\hline \hline
&&&&\\[-8pt]
\ ${\mbox{\tt
\fontsize{10pt}{0mm}\selectfont
E}}_6$&{\mbox{\tt
\fontsize{10pt}{0mm}\selectfont
EI}}=${\mbox{\tt
\fontsize{10pt}{0mm}\selectfont
E}}_{6(6)}$&
${\mathfrak {sp}}_4$&36&42\\
\ ${\mbox{\tt
\fontsize{10pt}{0mm}\selectfont
E}}_6$&{\mbox{\tt
\fontsize{10pt}{0mm}\selectfont
EII}}=${\mbox{\tt
\fontsize{10pt}{0mm}\selectfont
E}}_{6(2)}$&${\mathfrak {su}}_2
\oplus{\mathfrak {su}}_6$&38&40\\
\ ${\mbox{\tt
\fontsize{10pt}{0mm}\selectfont
E}}_6$&{\mbox{\tt
\fontsize{10pt}{0mm}\selectfont
EIII}}=${\mbox{\tt
\fontsize{10pt}{0mm}\selectfont
E}}_{6(-14)}$&${\mathfrak {so}}_{10}
\oplus{\mathbb R}$&46&32\\
\ ${\mbox{\tt
\fontsize{10pt}{0mm}\selectfont
E}}_6$&{\mbox{\tt
\fontsize{10pt}{0mm}\selectfont
EIV}}=${\mbox{\tt
\fontsize{10pt}{0mm}\selectfont
E}}_{6(-26)}$& ${\mbox{\tt
\fontsize{10pt}{0mm}\selectfont F}}_{4(-52)}$&52&26\\[3pt]
\hline
&&&&\\[-10pt]
\ ${\mbox{\tt
\fontsize{10pt}{0mm}\selectfont
E}}_7$&{\mbox{\tt
\fontsize{10pt}{0mm}\selectfont
EV}}=${\mbox{\tt
\fontsize{10pt}{0mm}\selectfont
E}}_{7(7)}$&
${\mathfrak {su}}_8$&63&70\\
\ ${\mbox{\tt
\fontsize{10pt}{0mm}\selectfont
E}}_7$&{\mbox{\tt
\fontsize{10pt}{0mm}\selectfont
EVI}}=${\mbox{\tt
\fontsize{10pt}{0mm}\selectfont
E}}_{7(-5)}$& ${\mathfrak {su}}_2
\oplus{\mathfrak {so}}_{12}$&69&64\\
\ ${\mbox{\tt
\fontsize{10pt}{0mm}\selectfont
E}}_7$&{\mbox{\tt
\fontsize{10pt}{0mm}\selectfont
EVII}}=${\mbox{\tt
\fontsize{10pt}{0mm}\selectfont
E}}_{7(-25)}$ &${\mbox{\tt
\fontsize{10pt}{0mm}\selectfont
E}}_{6(-78)}
\oplus{\mathbb R}$&79&54\\[3pt]
\hline
&&&&\\[-10pt]
\ ${\mbox{\tt
\fontsize{10pt}{0mm}\selectfont
E}}_8$&{\mbox{\tt
\fontsize{10pt}{0mm}\selectfont
EVIII}}=${\mbox{\tt
\fontsize{10pt}{0mm}\selectfont
E}}_{8(8)}$&
${\mathfrak {so}}_{16}$&120&128\\
\ ${\mbox{\tt
\fontsize{10pt}{0mm}\selectfont
E}}_8$&{\mbox{\tt
\fontsize{10pt}{0mm}\selectfont
EIX}}=${\mbox{\tt
\fontsize{10pt}{0mm}\selectfont
E}}_{8(-24)}$ &${\mathfrak
{su}}_2\oplus{\mbox{\tt
\fontsize{10pt}{0mm}\selectfont
E}}_{7(-133)}$
&136&112\\[3pt]
\hline
&&&&\\[-10pt]
\ ${\mbox{\tt
\fontsize{10pt}{0mm}\selectfont
F}}_4$&{\mbox{\tt
\fontsize{10pt}{0mm}\selectfont
FI}}=${\mbox{\tt
\fontsize{10pt}{0mm}\selectfont
F}}_{4(4)}$& ${\mathfrak {su}}_2
\oplus{\mathfrak {sp}}_3$&24&28\\
\ ${\mbox{\tt
\fontsize{10pt}{0mm}\selectfont
F}}_4$&{\mbox{\tt
\fontsize{10pt}{0mm}\selectfont
FII}}=${\mbox{\tt
\fontsize{10pt}{0mm}\selectfont
F}}_{4(-20)}$
&${\mathfrak {so}}_9$&36&16\\[3pt]
\hline
&&&&\\[-10pt]
\ ${\mbox{\tt
\fontsize{10pt}{0mm}\selectfont
G}}_2$&{\mbox{\tt
\fontsize{10pt}{0mm}\selectfont
GI}}=${\mbox{\tt
\fontsize{10pt}{0mm}\selectfont
G}}_{2(2)}$&${\mathfrak {so}}_3
\oplus{\mathfrak {so}}_3$&6&8\\
\end{tabular}
}
\end{center}

\vskip 3mm

%%%%%%%%%%%%%%%%%%%%%%%%%%%%%%%%%%%%%%%%
Recall the classical approach to
classifying nilpotent orbits by means of
their characte\-ristics and weighted
Dynkin diagrams, cf.\,\cite{CM}, \cite{M}.
Fix a $\theta$-stable Cartan subalgebra
$\te$ of $\g$ such that $\te\cap \ka$ is a
Cartan subalgebra of $\ka$. Denote
$$
l:=\dim^{}_{\C} \te, \quad s:=\dim^{}_{\C}
(\te\cap \ka).
$$
Let $\Delta^{}_{(\g, \te)}$ and
$\Delta^{}_{(\ka, \te\cap \ka)}$ be the
root systems of resp.\;\hskip -.4mm $(\g,
\te)$ and $(\ka, \te\cap \ka)$. Given a
nonzero element $e\in \cNg$
(resp.,\;\hskip -.4mm\,$e\in \cNp$), there
are the elements $h, f\in \g$
(resp.,\;\hskip -.4mm\,$h\in \ka$, $f\in
\pe$) such that $\{e, h, f\}$ is an
${\mathfrak {sl}}_2$-triple. The
intersection of $\te$ (resp.,\;\hskip
-.4mm\,$\te\cap \ka$) with the orbit
${\cal O}=G\cdot e$ (resp.,\;\hskip
-.4mm\,$K\cdot e$) contains a unique
element $h_0$ lying in a fixed Weyl
chamber of $\Delta^{}_{(\g, \te)}$
(resp.,\;\hskip -.4mm\,$\Delta^{}_{(\ka,
\te\cap \ka)}$). The mapping ${\cal
O}\mapsto h_0$ is a well defined injection
of the set of all nonzero $G$-orbits in
$\cNg$ (resp.,\;\hskip -.4mm\,$K$-orbits
in $\cNp$) into $\te$ (resp.,\;\hskip
-.4mm\,$\te\cap \ka$). The image $h_0$ of
the orbit ${\cal O}$ under this injection
is called the {\it characteristic} of
$\cal O$ (and $e$). Thus the $G$-orbits in
$\cNg$ (resp.,\;\hskip -.4mm\;$K$-orbits
in $\cNp$) are defined by their
characteristics. In turn, the
characteristics are defined by the
numerical data, namely, the system of
values $\{\beta_j(h_0)\}$ where
$\{\beta_j\}$ is a fixed basis of $\te^*$
(resp.,\;\hskip -.4mm\,$(\te\cap \ka)^*$).
In practice, $\{\beta_j\}$ is always
chosen to be a base of some root system.
Then assigning the integer $\beta_j(h_0)$
to each node $\beta_j$ of the Dynkin
diagram of this base yields the {\it
weighted Dynkin diagram} of the orbit
${\cal O}$, denoted by ${\rm Dyn}\,{\cal
O}$. It uniquely defines ${\cal O}$. This
choice of the base $\{\beta_j\}$ is the
following.

For $\te^*$, since $\g$ is semisimple, it
is natural to take $\{\beta_j\}$ to be the
base $\alpha_1,\ldots,\alpha_l$ of
$\Delta^{}_{(\g, \te)}$ defining the fixed
Weyl chamber that is used in the
definition of characteristics.  Below the
Dynkin diagrams of $G$-orbits in $\cNg$
are considered with respect to this base
$\alpha_1,\ldots,\alpha_l$.

To describe $\{\beta_j\}$ for $(\te\cap
\ka)^*$, denote by $\alpha_0$ the lowest
root of $\Delta^{}_{(\g, \te)}$. The
extended Dynkin diagrams of
$\Delta^{}_{(\g, \te)}$ indicating the
location of the $\alpha_i$'s are given in
Table 6.

The real form $\gR$ of $\g$ is of inner
type (i.e., $s=l$, $\te\subset \ka$, and
thereby $\Delta_{(\ka, \te\cap
\ka)}\subset \Delta_{(\g, \te)}$) if and
only if $\gR\neq {\mbox{\tt
\fontsize{13pt}{0mm}\selectfont
E}}_{6(6)}$, ${\mbox{\tt
\fontsize{13pt}{0mm}\selectfont
E}}_{6(-26)}$. If $\gR$ is of inner type,
then $\ka$ is semisimple if and only if
$\gR\neq {\mbox{\tt
\fontsize{13pt}{0mm}\selectfont
E}}_{6(-14)}, {\mbox{\tt
\fontsize{13pt}{0mm}\selectfont
E}}_{7(-25)}$; in the last two cases the
center of $\ka$ is one-dimensional. If
$\gR$ is of inner type, then the set
$\{\alpha_0,\ldots,\alpha_l\}$ contains a
unique base $\beta_1,\ldots,\beta_r$ of
$\Delta_{(\ka, \te\cap \ka)}$. Moreover,
$\alpha_0$ is one of the $\beta_i$'s
$\Longleftrightarrow$ $\ka$ is semisimple
$\Longleftrightarrow$ $r=l$. If $\gR$ is
of inner type and $\ka$ is not semisimple,
then $r=l-1$ and
$\beta_j\in\{\alpha_1,\ldots,\alpha_l\}$
for all $j=1,\ldots, r$; in this case we
set $\beta_l:=\alpha_l$.

Thus in all cases where $\gR$ is of inner
type, we have a basis $\{\beta_j\}$ of
$(\te\cap\ka)^*$ that is the system of
simple roots of some root system. The
location of the $\beta_j$'s is given in
Table~6.

%%\

%%\vskip -9mm

%%\

\begin{center}
\centerline{\text TABLE
6}\nopagebreak\vskip -1mm \centerline{ }
\leavevmode \epsfxsize =15cm
%%15cm
\epsffile{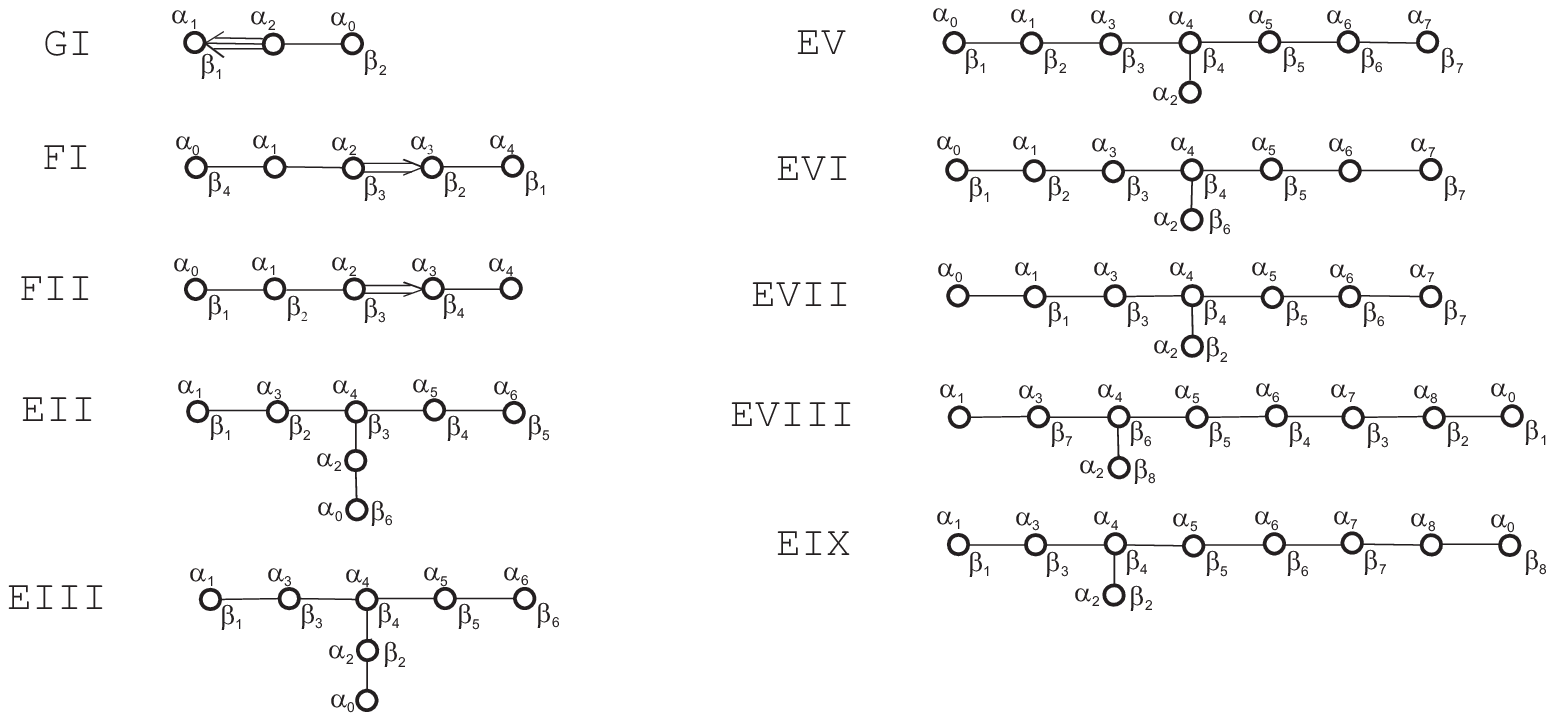}
%%\centerline{\text{ ¨á.\,1}}
\end{center}

If $\gR$ is of outer type, then $\ka$ is
semisimple, so we may, and will, take
$\{\beta_j\}$ to be a base of
$\Delta_{(\ka, \te\cap \ka)}$. We take the
following base. If $\gR={\mbox{\tt
\fontsize{13pt}{0mm}\selectfont
E}}_{6(-26)}$, then
$$
\beta_4=\alpha_1|_{\te\cap
\ka}=\alpha_6|_{\te\cap \ka},\hskip 4mm
\beta_3=\alpha_3|_{\te\cap \ka}
=\alpha_5|_{\te\cap \ka},\hskip 4mm
\beta_2=\alpha_4|_{\te\cap \ka},\hskip 4mm
\beta_1=\alpha_2|_{\te\cap \ka},
$$
with the Dynkin diagram

\begin{center}
\leavevmode \epsfxsize =3cm
\epsffile{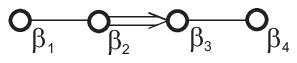}
\end{center}

\noindent and if $\gR={\mbox{\tt
\fontsize{13pt}{0mm}\selectfont
E}}_{6(6)}$, then
\begin{gather*}
\begin{gathered}
\beta_1=-2\beta_2-3\beta_3-2\beta_4-
\alpha_2|_{\te\cap \ka},\\
 \beta_2=\alpha_1|_{\te\cap
\ka}=\alpha_6|_{\te\cap \ka},\hskip 4mm
\beta_3=\alpha_3|_{\te\cap \ka}
=\alpha_5|_{\te\cap \ka},\hskip 4mm
\beta_4=\alpha_4|_{\te\cap \ka}
\end{gathered}
\end{gather*}
with the Dynkin diagram

\begin{center}
\leavevmode \epsfxsize =3cm
\epsffile{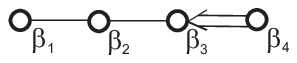}
\end{center}

Below the weighted Dynkin diagrams of
$K$-orbits in $\cNg$ are considered with
respect to the described base
$\{\beta_j\}$.

\medskip

In \cite{D3}, \cite{D4}  one finds:

\smallskip

\blist

\item[\rm (D1)] the weighted Dynkin
diagrams of all nonzero orbits $G\cdot x$
and $ K\cdot x$ for $x\!\in\! {\cal
N}(\pe)$,

\item[\rm (D2)] the type and dimension of
the reductive Levi factor of $\z_{\ka}(x)$
and $\dim^{}_{\C}\z_{\ka}(x)$,

\item[\rm (D3)] the type of the reductive
Levi factor of ${\z}_{\gR}(y)$ for an
element $y$ of the $\GRc$-orbit in $\cNgr$
corresponding to $K\cdot x$ via the
Kostant--Sekiguchi bijection.

\end{list}

\smallskip

Given the weighted Dynkin diagram of
$G\cdot x$, one finds the type of
reductive Levi factor of $\z_{\g}(x)$ in
\cite{El}, \cite{Ca}. So using (D1), one
can apply Theorem \ref{levi} for explicit
classifying $(-1)$-distinguished
$K$-orbits in $\cNp$. On the other hand,
given (D3), Theorem \ref{kost-sek} can be
applied for this purpose as well. So
following either of these ways, one
obtains, for every exceptional simple
$\g$, the explicit classification of
$(-1)$-distinguished $K$-orbits in $\cNp$
in terms of their weighted Dynkin
diagrams. The final result is the
following.

\begin{theorem}\label{except}
For all exceptional simple Lie algebras
$\g$ and all conjugacy classes of elements
$\theta\in {\rm Aut}\,\g$ of order~$2$,
all $(-1)$-distinguished $K$-orbits in
$\cNp$ are listed in Tables $7$--$18$ at
the end of this paper.
%%for
%%all exceptional simple Lie algebras $\g$
%%and all conjugacy classes of elements
%%$\theta\in {\rm Aut}\,\g$ of order~$2$.
\end{theorem}

Tables 7--18 contain further information
as summarized below.

\smallskip

\blist

\item[\rm (T1)] The conjugacy class of
$\theta$ is defined by specifying the type
of noncompact real form $\gR$ of $\g$
canonically corresponding to this class.

\item[\rm (T2)] Column 2 gives the weights
$\{\beta_j(h)\}$ (listed in the order of
increasing of $j$) of the weighted Dynkin
diagram ${\rm Dyn}\,K\cdot x$ of the
$(-1)$-distinguished $K$-orbit $K\cdot x$
in $\cNp$, where $h$ is the characteristic
of $K\cdot x$.

\item[\rm (T3)] Column 3 gives the weights
$\{\alpha_j(H)\}$ (listed in the order of
increasing of $j$) of the weighted Dynkin
diagram ${\rm Dyn}\,G\cdot x$ of the
$G$-orbit $G\cdot x$ in $\cNg$, where $H$
is the characteristic of $G\cdot x$.

\item[\rm (T4)] Column 4 gives
$\dim^{}_{\C} K\cdot x$, and hence, by
\eqref{2}, also $\dim^{}_{\C} G\cdot x$.

\item[\rm (T5)] Column 5 gives the number
of $K$-orbits (not necessarily
$(-1)$-distinguished) in $G\cdot x\cap
\pe$. One of them is $K\cdot x$.
  Since this number is equal
  to the number of
$K$-orbits $K\cdot x'$ in $\cNp$ such that
${\rm Dyn}\,G\cdot x'={\rm Dyn}\,G\cdot
x$, one finds it using the tab\-les in
\cite{D3}, \cite{D4}.

\item[\rm (T6)] Columns 6 gives the type
of the reductive Levi factor  of
$\z_{\ka}(x)$.
 By ${\mbox{\tt \fontsize{13pt}{0mm}\selectfont T}}_m$ is denoted the Lie algebra of
an $m$-dimensional torus.

\item[\rm (T7)] Column 7 gives the complex
dimension of the unipotent radical of
$\z_{\ka}(x)$.

\end{list}

\smallskip

\noindent{\it Remark 1.} For $\mbox{\tt
\fontsize{13.1pt}{0mm}\selectfont
E}_{6(6)}$ and $\mbox{\tt
\fontsize{13.1pt}{0mm}\selectfont
E}_{6(-26)}$, our numeration of
$\{\alpha_j\}$ is the same as in
\cite{D3}; it differs from that in
\cite{D4}. For $\mbox{\tt
\fontsize{13.1pt}{0mm}\selectfont
E}_{6(-26)}$, our numeration of
$\{\beta_j\}$ differs from that in
\cite{D3}.

\medskip

\noindent{\it Remark 2.} One finds   in
\cite{D5}, \cite{D6}, \cite{D7} the
explicit classification of all Cayley
triples in $\g$. This yields the explicit
representatives of all $K$-orbits in
$\cNp$ and $\GRc$-orbits in $\cNgr$.

\section*{6. Geometric properties
 \label{exept}}

Geometric properties of the varieties
$\mathbf P(\ov{\cal O})$ where ${\cal O}$
is a nilpotent $K$-orbit in $\cNp$ have
been studies by several authors. Thereby
their results provide some information on
the geometry of the projective self-dual
varieties that we associated with
symmetric spaces. However these results
are less complete than in the adjoint case
(where a
%%see \cite{P2} for a
rather detailed information about singular
loci and normality of the projective
self-dual varieties $\mathbf P(\ov{\cal
O})$ is available, see \cite{P2}). Some
phenomena valid for the projective
self-dual $\mathbf P(\ov{\cal O})$'s in
the adjoint case fail in general. Below we
briefly summarize some facts about the
geometry of the projective self-dual
varieties that we associated with
symmetric spaces.

\medskip

{\it Intersections of $G$-orbits in
${\mathbf P}(\cNg)$ with the linear
subspace ${\mathbf P}(\pe)$ }

\smallskip

Let $\mathcal O$ be a nonzero  $G$-orbit
in $\cNg$, and $X={\mathbf
P}(\overline{\mathcal O})$. If $X\cap
{\mathbf P}(\pe)\neq \varnothing$, then,
%%the subspace ${\mathbf P}(\pe)$ of
%%${\mathbf P}(\g)$ intersects $X$, then,
by \eqref{2}, all irreducible components
$Y_1,\ldots Y_s$ of the variety $X\cap
{\mathbf P}(\pe)$
 have dimension $\frac{1}{2}\dim X$, and
 every $Y_j$ is the closure of a
 $K$-orbit in ${\mathbf P}(\cNp)$.
Theorems \ref{dist}, \ref{-dist} and
Definition \ref{-1} imply that
$$ X \mbox{ is
self-dual }\hskip 1.5mm
\Longrightarrow\hskip 1mm \mbox{ all }
Y_1, \ldots, Y_s \mbox{ are self-dual. }
$$
There are many examples showing that the
converse in not true. For instance, one
deduces from Table 8 that $\g={\mbox{\tt
\fontsize{12pt}{0mm}\selectfont E}_6}$,
$\g_{\R}={\mbox{\tt
\fontsize{12pt}{0mm}\selectfont
E}_{6(2)}}$ and ${\rm Dyn}\hskip
.5mm{\mathcal O}=111011, 121011$ or
$220202$ are such cases. There are also
many examples where some of the $Y_j$'s
are self-dual and some are not: for
instance, this is so if $\g={\mbox{\tt
\fontsize{12pt}{0mm}\selectfont E}_6}$,
$\g_{\R}={\mbox{\tt
\fontsize{12pt}{0mm}\selectfont
E}_{6(2)}}$ and ${\rm Dyn}\hskip
.5mm{\mathcal O}=020000, 001010, 000200,
020200$ or $220002$. Finally, there are
many instances where all $Y_j$ are not
self-dual (to obtain them,
e.g.,\;compare\;Tables 7--18 with tables
in \cite{D3}, \cite{D4}).

\medskip

{\it Affine $(-1)$-distinguished orbits}

\smallskip

Recall that by Matsushima's criterion, an
orbit of a reductive algebraic group
acting on an affine algebraic variety is
affine if and only if the stabilizer of a
point of this orbit is reductive,
cf.,\;e.g.\,\cite[Theorem 4.17]{PV}.
Therefore distinguished $G$-orbits in
$\cNg$ are never affine. However in
general there exist affine
$(-1)$-distinguished $K$-orbits in $\cNp$.
For instance, if $\g$ is exceptional
simple, we immediately obtain their
classification from Tables 7--18: these
are precisely the orbits for which the
number in the last column is~$0$.

If $K\cdot x$ is affine, then each
irreducible component of the boundary
${\mathbf P}(\overline{K\cdot x})\setminus
{\mathbf P}(K\cdot x)$ has codimension 1
in ${\mathbf P}(\overline{K\cdot x})$,
cf.\,\cite[Lemma 3]{P1}. Therefore if a
point of the open $K$-orbit of such an
irreducible component lies in the singular
locus of ${\mathbf P}(\overline{K\cdot
x})$, then ${\mathbf P}(\overline{K\cdot
x})$ is not normal.

\medskip

{\it Orbit closure ordering and orbit
decomposition of the orbit boundary
${\mathbf P}(\overline{K\cdot x})\setminus
{\mathbf P}(K\cdot x)$}

\smallskip

The closure ordering on the set of
$K$-orbits in $\cNp$ (resp.,\;\hskip
-.4mm\;${\rm Ad}(\g^{}_{\R})$-orbits in
$\cNgr$) is defined by the condition that
$ {\cal O}_1>{\cal O}_2$ if and only if
${\cal O}_2$ is contained in the closure
of ${\cal O}_1$ and ${\cal O}_1\neq {\cal
O}_2$. According to \cite{BS}, the
Kostant--Sekiguchi bijection preserves the
closure ordering. Clearly its describing
can be reduced to the case of simple $\g$.
In this case, the explicit description of
the closure ordering is obtained in
\cite{D1}, \cite{D2}, \cite{D8},
\cite{D9}, \cite{D10}, \cite{D11},
\cite{D12}, \cite{D13} (see also \cite{O},
\cite{Se}).

This information and Theorems
\ref{sl-compa}--\ref{except} yield the
orbit decomposition of every irreducible
component of the orbit boundary ${\mathbf
P}(\overline{K\cdot x})\setminus {\mathbf
P}(K\cdot x)$ for any $(-1)$-distinguished
$K$-orbit $K\cdot\nobreak x$. In
particular, this gives the dimensions of
these components and their intersection
configurations.

\medskip

{\it Singular locus of a self-dual variety
${\mathbf P}(\overline {K\cdot x})$}

\smallskip

Apart from a rather detailed information
on the singular loci of the self-dual
projectivized nilpotent orbit closures in
the adjoint case (see \cite{P2}), in
general case only a partial information is
available. In \cite{O}, for $(\g,
\ka)=(\ssl_n(\C), \so_n(\C))$ and
$(\ssl_{2n}(\C), \spl_n(\C))$, the
nor\-ma\-li\-ty of the nil\-po\-tent orbit
closures $\overline{K\cdot x}$ (hence that
of ${\mathbf P}(\overline{K\cdot x})$) is
studied. In par\-ti\-cu\-lar, it is shown
there that for $(\g, \ka)=(\ssl_{2n}(\C),
\spl_n(\C))$, normality holds for any
$\overline{K\cdot x}$, but if $(\g,
\ka)=(\ssl_n(\C), \so_n(\C))$, this is not
the case. See also \cite{Se}, \cite{O},
\cite{SeSh}, where the local equations of
the generic singularities of some orbit
closures in ${\mathbf P}(\cNp)$ are found.
By \cite[Proposition 4]{P3}, Hesselink's
desingularization of the closures of
Hesselink strata, \cite{He},
cf.\,\cite{PV}, yields a desingularization
of any orbit closure ${\mathbf
P}(\overline{K\cdot x})$ in ${\mathbf
P}(\cNp)$ (in \cite{R}, another approach
to desingularization of such orbit
closured is considered).

\pagebreak

{\fontsize{9pt}{6mm}\selectfont
\begin{center}{TABLE 7} \\
{ $(-1)$-distinguished $K$-orbits $
%%{\cal O}=
K\cdot x$ in $\cNp$ for $\gR=
{\mbox{\tt \fontsize{10pt}{0mm}\selectfont E}}_{6(6)}$}\\
\
\\
\begin{tabular}{c|c|c|c|c|c|c}

No.&${\rm Dyn}\,
%%{\cal O}
K\cdot x$ &${\rm Dyn}\, G\cdot x$ &
$\dim^{}_{\C} K\cdot x$
%% {\cal O}$
& $\sharp (G\cdot x\cap \pe)/K$ &
$\z_{\ka}(x)/{\rm rad}_u\z_{\ka}(x)$&
$\dim^{}_{\C}{\rm rad}_u\z_{\ka}(x)$
\\
\hline \hline
&&&&&&\\[-13pt]
1&2222&202222 &35&1&0&1\\
\hline
2&2202&220002 &33&2&0&3\\[-3pt]
3&0220&220002 &33&2&0&3\\
\hline
4&4224&222222 &36&1&0&0\\
\end{tabular}
\end{center}
}

\vskip 6mm
%%\newpage
%%\pagebreak

 {\fontsize{9pt}{6mm}\selectfont
\begin{center}{TABLE 8} \\
{$(-1)$-distinguished $K$-orbits $
%%{\cal O}=
K\cdot x$ in $\cNp$ for $\gR={\mbox{\tt
\fontsize{10pt}{0mm}\selectfont E}}_{6(2)}$}\\
\
\\
\begin{tabular}{c|c|c|c|c|c|c}

No. &${\rm Dyn}\,K\cdot x$
%%{\cal O}$
&${\rm Dyn}\,G\cdot x$

& $\dim^{}_{\C} K\cdot x$
%%{\cal O}$
& $\sharp (G\cdot x\cap \p)/K$ &
$\z_{\ka}(x)/{\rm rad}_u\z_{\ka}(x)$&
$\dim^{}_{\C}{\rm rad}_u\z_{\ka}(x)$
\\
\hline \hline
&&&&&& \\[-13pt]
1&000004&020000&21&3&$2{\mbox{\tt
\fontsize{10pt}{0mm}\selectfont A}}_2$& 1
\\
\hline 2&301000&001010&25&3&${\mbox{\tt
\fontsize{10pt}{0mm}\selectfont
A}}_1+{\mbox{\tt
\fontsize{10pt}{0mm}\selectfont T}}_1$&9
\\[-3pt]
3&001030&001010&25&3&${\mbox{\tt
\fontsize{10pt}{0mm}\selectfont
A}}_1+{\mbox{\tt
\fontsize{10pt}{0mm}\selectfont T}}_1$&9
\\
\hline 4&004000&000200&29&3&${\mbox{\tt
\fontsize{10pt}{0mm}\selectfont T}}_2$&7\\[-3pt]
5&020204&000200&29&3&${\mbox{\tt
\fontsize{10pt}{0mm}\selectfont T}}_2$&7\\
\hline 6&004008&020200&30&2&${\mbox{\tt
\fontsize{10pt}{0mm}\selectfont A}}_2$&0
\\
\hline 7&400044&220002&30&2&${\mbox{\tt
\fontsize{10pt}{0mm}\selectfont
A}}_1+{\mbox{\tt
\fontsize{10pt}{0mm}\selectfont T}}_1$&4
\\
\hline 8&121131&111011&31&2&${\mbox{\tt
\fontsize{10pt}{0mm}\selectfont T}}_1$&6\\[-3pt]
9&311211&111011&31&2&${\mbox{\tt
\fontsize{10pt}{0mm}\selectfont T}}_1$&6\\
\hline 10&313104&121011&32&2&${\mbox{\tt
\fontsize{10pt}{0mm}\selectfont T}}_1$&5\\[-3pt]
11&013134&121011&32&2&${\mbox{\tt
\fontsize{10pt}{0mm}\selectfont T}}_1$&5\\
\hline
12&222222&200202&33&2&0&5\\[-3pt]
13&040404&200202&33&2&0&5\\ \hline
14&224224&220202&34&2&${\mbox{\tt
\fontsize{10pt}{0mm}\selectfont T}}_1$&3\\[-3pt]
15&404048&220202&34&2&${\mbox{\tt
\fontsize{10pt}{0mm}\selectfont T}}_1$&3\\
\hline 16&440444&222022&35&1&0&3\\ \hline
17&444448&222222&36&1&0&2\\[-3pt]
\end{tabular}
\end{center}
}

\vskip 5mm

{\fontsize{9pt}{6mm}\selectfont
\begin{center}{TABLE 9} \\
{ $(-1)$-distinguished $K$-orbits $
%%{\cal O}=
K\cdot x$ in $\cNp$ for $\gR={\mbox{\tt
\fontsize{10pt}{0mm}\selectfont E}}_{6(-26)}$}\\
\
\\
\begin{tabular}{c|c|c|c|c|c|c}

No.&${\rm Dyn}\,K\cdot x$
%%{\cal O}$
&${\rm Dyn}\, G\cdot x$

& $\dim^{}_{\C} K\cdot x$
%%{\cal O}$
& $\sharp (G\cdot x\cap \p)/K$ &
$\z_{\ka}(x)/{\rm rad}_u\z_{\ka}(x)$&
$\dim^{}_{\C}{\rm rad}_u\z_{\ka}(x)$
\\
\hline \hline
&&&&&&\\[-13pt]
1&  0001&100001 &16&1&${\mbox{\tt
\fontsize{10pt}{0mm}\selectfont B}}_3$&15
\\
\hline 2&  0002&200002 &24&1&${\mbox{\tt
\fontsize{10pt}{0mm}\selectfont G}}_2$&14
\\
\end{tabular}
\end{center}
}

\vskip 5mm

{\fontsize{9pt}{6mm}\selectfont
\begin{center}{TABLE } 10\\
{ $(-1)$-distinguished $K$-orbits $
%%{\cal O}=
K\cdot x$ in $\cNp$ for $\gR={\mbox{\tt
\fontsize{10pt}{0mm}\selectfont E}}_{6(-14)}$}\\
\
\\
\begin{tabular}{c|l|c|c|c|c|c}

No.&\ ${\rm Dyn}\,K\cdot x$
%%{\cal O}$
&${\rm Dyn}\,G\cdot x$

& $\dim^{}_{\C}K\cdot x$
%%{\cal O}$
& $\sharp (G\cdot x\cap \p)/K$ &
$\z_{\ka}(x)/{\rm rad}_u\z_{\ka}(x)$&
$\dim^{}_{\C}{\rm rad}_u\z_{\ka}(x)$
\\
\hline \hline
&&&&&&\\[-13pt]
1&100001&100001&16&3&${\mbox{\tt
\fontsize{10pt}{0mm}\selectfont
B}}_3+{\mbox{\tt
\fontsize{10pt}{0mm}\selectfont T}}_1$& 8
  \\[-3pt]
2&10000$-2$&100001&16&3&${\mbox{\tt
\fontsize{10pt}{0mm}\selectfont
B}}_3+{\mbox{\tt
\fontsize{10pt}{0mm}\selectfont T}}_1$& 8
\\
\hline 3&10101$-2$&110001&23&2&${\mbox{\tt
\fontsize{10pt}{0mm}\selectfont
A}}_2+{\mbox{\tt
\fontsize{10pt}{0mm}\selectfont T}}_1$& 14
\\[-3pt]
4&11100$-3$&110001&23&2&${\mbox{\tt
\fontsize{10pt}{0mm}\selectfont
A}}_2+{\mbox{\tt
\fontsize{10pt}{0mm}\selectfont T}}_1$& 14
\\
\hline 5&40000$-2$&200002&24&1&${\mbox{\tt
\fontsize{10pt}{0mm}\selectfont G}}_2$&8
\\
\hline 6&03001$-2$&120001&26&2&${\mbox{\tt
\fontsize{10pt}{0mm}\selectfont
B}}_2+{\mbox{\tt
\fontsize{10pt}{0mm}\selectfont T}}_1$&9
\\[-3pt]
7&01003$-6$&120001&26&2&${\mbox{\tt
\fontsize{10pt}{0mm}\selectfont
B}}_2+{\mbox{\tt
\fontsize{10pt}{0mm}\selectfont T}}_1$&9
\\
\hline 8&02202$-6$&220002&30&1&${\mbox{\tt
\fontsize{10pt}{0mm}\selectfont
A}}_1+{\mbox{\tt
\fontsize{10pt}{0mm}\selectfont T}}_1$&12
\\
\end{tabular}
\end{center}
}

\vskip 8mm

{\fontsize{9pt}{6mm}\selectfont
\begin{center}{TABLE 11} \\
{ $(-1)$-distinguished $K$-orbits $
%%{\cal O}=
K\cdot x$ in $\cNp$ for $\gR={\mbox{\tt
\fontsize{10pt}{0mm}\selectfont E}}_{7(-5)}$}\\
\
\\
\begin{tabular}{c|c|c|c|c|c|c}
No.&${\rm Dyn}\,K\cdot x$
%%{\cal O}$
& ${\rm Dyn}\,G\cdot x$&
$\dim^{}_{\C}K\cdot x$
%%{\cal O}$
& $\sharp (G\cdot x\cap \p)/K$ &
$\z_{\ka}(x)/{\rm rad}_u\z_{\ka}(x)$&
$\dim^{}_{\C}{\rm rad}_u\z_{\ka}(x)$
\\
\hline \hline
&&&&&&\\[-13pt]
1&0000004&2000000&33&3&${\mbox{\tt
\fontsize{10pt}{0mm}\selectfont A}}_5$&1
\\
\hline 2&4000000&0000020&42&2&${\mbox{\tt
\fontsize{10pt}{0mm}\selectfont
G}}_2+{\mbox{\tt
\fontsize{10pt}{0mm}\selectfont A}}_1$&10
\\
\hline 3&0000400&0020000&47&3&$3{\mbox{\tt
\fontsize{10pt}{0mm}\selectfont A}}_1$&13
\\[-3pt]
4&0002004&0020000&47&3&$3{\mbox{\tt
\fontsize{10pt}{0mm}\selectfont A}}_1$&13
\\
\hline 5&0000408&2020000&48&2&${\mbox{\tt
\fontsize{10pt}{0mm}\selectfont C}}_3$&0
\\
\hline 6&2010112&0001010&49&1&${\mbox{\tt
\fontsize{10pt}{0mm}\selectfont
A}}_1+{\mbox{\tt
\fontsize{10pt}{0mm}\selectfont T}}_1$&16
\\
\hline 7&0400004&2000020&50&2&${\mbox{\tt
\fontsize{10pt}{0mm}\selectfont
A}}_2+{\mbox{\tt
\fontsize{10pt}{0mm}\selectfont T}}_1$&10
\\
\hline 8&1111101&1001010&52&1&${\mbox{\tt
\fontsize{10pt}{0mm}\selectfont T}}_2$&15  \\
\hline 9&2010314&2001010&53&1&${\mbox{\tt
\fontsize{10pt}{0mm}\selectfont
A}}_1+{\mbox{\tt
\fontsize{10pt}{0mm}\selectfont T}}_1$&12
\\
\hline 10&0040000&0002000&53&1&${\mbox{\tt
\fontsize{10pt}{0mm}\selectfont A}}_1$&13
\\
\hline 11&0202202&0020020&55&2&${\mbox{\tt
\fontsize{10pt}{0mm}\selectfont A}}_1$&11
\\[-3pt]
12&0004004&0020020&55&2&${\mbox{\tt
\fontsize{10pt}{0mm}\selectfont A}}_1$&11
\\
\hline
13&0202404&2020020&56&2&$2{\mbox{\tt
\fontsize{10pt}{0mm}\selectfont A}}_1$&7
\\[-3pt]
14&0400408&2020020&56&2&$2{\mbox{\tt
\fontsize{10pt}{0mm}\selectfont A}}_1$&7
\\
\hline 15&4004000&0002020&57&1&${\mbox{\tt
\fontsize{10pt}{0mm}\selectfont A}}_1$&9
\\
\hline 16&0404004&2002020&59&1&${\mbox{\tt
\fontsize{10pt}{0mm}\selectfont T}}_1$&9  \\
\hline 17&0404408&2022020&60&1&${\mbox{\tt
\fontsize{10pt}{0mm}\selectfont A}}_1$&6
\\
\end{tabular}
\end{center}
}

\newpage

\vskip 8mm

{\fontsize{9pt}{6mm}\selectfont
\begin{center}{TABLE 12} \\
{$(-1)$-distinguished $K$-orbits $
%%{\cal O}=
K\cdot x$ in $\cNp$ for $\gR={\mbox{\tt
\fontsize{10pt}{0mm}\selectfont E}}_{7(7)}$}\\
\
\\
\begin{tabular}{c|c|c|c|c|c|c}

No. &${\rm Dyn}\,K\cdot x$
%%{\cal O}$
&${\rm Dyn}\,G\cdot x$ & $\dim^{}_{\C}
K\cdot x$
%%{\cal O}$
& $\sharp (G\cdot x\cap \p)/K$ &
$\z_{\ka}(x)/{\rm rad}_u\z_{\ka}(x)$&
$\dim^{}_{\C}{\rm rad}_u\z_{\ka}(x)$
\\
\hline \hline
&&&&&&\\[-13pt]
1&4000000&0200000&42&4&${\mbox{\tt
\fontsize{10pt}{0mm}\selectfont G}}_2$&7
\\[-3pt]
2&0000004&0200000&42&4&${\mbox{\tt
\fontsize{10pt}{0mm}\selectfont G}}_2$&7
\\
\hline 3&0040000&0000200&50&4&${\mbox{\tt
\fontsize{10pt}{0mm}\selectfont A}}_1$&10
\\[-3pt]
4&0000400&0000200&50&4&${\mbox{\tt
\fontsize{10pt}{0mm}\selectfont A}}_1$&10
\\
\hline 5&3101021&1001010&52&3&${\mbox{\tt
\fontsize{10pt}{0mm}\selectfont T}}_2$&9  \\[-3pt]
6&1201013&1001010&52&3&${\mbox{\tt
\fontsize{10pt}{0mm}\selectfont T}}_2$&9  \\
\hline 7&4004000&2000200&54&4&${\mbox{\tt
\fontsize{10pt}{0mm}\selectfont A}}_1$&6
\\[-3pt]
8&0004004&2000200&54&4&${\mbox{\tt
\fontsize{10pt}{0mm}\selectfont A}}_1$&6
\\
\hline
9&2220202&0002002&56&4&0&7  \\[-3pt]
10&2020222&0002002&56&4&0&7  \\[-3pt]
11&0400400&0002002&56&4&0&7  \\[-3pt]
12&0040040&0002002&56&4&0&7  \\
\hline
13&2222202&2002002&58&4&0&5  \\[-3pt]
14&2022222&2002002&58&4&0&5  \\[-3pt]
15&4004040&2002002&58&4&0&5  \\[-3pt]
16&0404004&2002002&58&4&0&5  \\
\hline 17&4220224&2002020&59&2&${\mbox{\tt
\fontsize{10pt}{0mm}\selectfont T}}_1$&3\\
\hline
18&2422222&2002022&60&4&0&3  \\[-3pt]
19&2222242&2002022&60&4&0&3  \\[-3pt]
20&4404040&2002022&60&4&0&3  \\[-3pt]
21&0404044&2002022&60&4&0&3  \\
\hline
22&4404404&2220202&61&2&0&2  \\[-3pt]
23&4044044&2220202&61&2&0&2  \\
\hline
24&4444044&2220222&62&2&0&1  \\[-3pt]
25&4404444&2220222&62&2&0&1  \\
\hline
26&8444444&2222222&63&2&0&0  \\[-3pt]
27&4444448&2222222&63&2&0&0  \\
\end{tabular}
\end{center}
}

\newpage

\vskip 5mm

{\fontsize{9pt}{6mm}\selectfont
\begin{center}{TABLE 13} \\
{$(-1)$-distinguished $K$-orbits $
%%{\cal O}=
K\cdot x$ in $\cNp$ for $\gR={\mbox{\tt
\fontsize{10pt}{0mm}\selectfont E}}_{7(-25)}$}\\
\
\\
\begin{tabular}{c|l|c|c|c|c|c}

No. &\ ${\rm Dyn}\,K\cdot x$
%%{\cal O}$
&${\rm Dyn}\, G\cdot x$

& $\dim^{}_{\C} K\cdot x$
%%{\cal O}$
& $\sharp (G\cdot x\cap \p)/K$ &
$\z_{\ka}(x)/{\rm rad}_u\z_{\ka}(x)$&
$\dim^{}_{\C}{\rm rad}_u\z_{\ka}(x)$
\\
\hline \hline
&&&&&&\\[-13pt]
1&0000002&0000002&27&4& ${\mbox{\tt
\fontsize{10pt}{0mm}\selectfont F}}_4$&0
\\[-3pt]
2&000000$-2$&0000002&27&4& ${\mbox{\tt
\fontsize{10pt}{0mm}\selectfont F}}_4$&0
\\
\hline
3&010010$-2$&1000010&38&2&${\mbox{\tt
\fontsize{10pt}{0mm}\selectfont
A}}_3+{\mbox{\tt
\fontsize{10pt}{0mm}\selectfont T}}_1$&25
\\[-3pt]
4&011000$-3$&1000010&38&2&${\mbox{\tt
\fontsize{10pt}{0mm}\selectfont
A}}_3+{\mbox{\tt
\fontsize{10pt}{0mm}\selectfont T}}_1$&25
\\
\hline
5&200002$-2$&2000002&43&4&${\mbox{\tt
\fontsize{10pt}{0mm}\selectfont B}}_3$&15
\\[-3pt]
6&400000$-2$&2000002&43&4&${\mbox{\tt
\fontsize{10pt}{0mm}\selectfont B}}_3$&15
\\[-3pt]
7&000004$-6$&2000002&43&4&${\mbox{\tt
\fontsize{10pt}{0mm}\selectfont B}}_3$&15
\\[-3pt]
8&200002$-6$&2000002&43&4&${\mbox{\tt
\fontsize{10pt}{0mm}\selectfont B}}_3$&15
\\
\hline
9&220002$-6$&2000020&50&1&${\mbox{\tt
\fontsize{10pt}{0mm}\selectfont
A}}_2+{\mbox{\tt
\fontsize{10pt}{0mm}\selectfont T}}_1$&20
\\
\hline
10&400004$-6$&2000022&51&2&${\mbox{\tt
\fontsize{10pt}{0mm}\selectfont G}}_2$&14
\\[-3pt]
11&400004$-10$&2000022&51&2&${\mbox{\tt
\fontsize{10pt}{0mm}\selectfont G}}_2$&14
\\
\end{tabular}
\end{center}
}

\vskip 5mm

{\fontsize{9pt}{6mm}\selectfont
\begin{center}{TABLE 14} \\
{$(-1)$-distinguished $K$-orbits $
%%{\cal O}=
K\cdot x$ in $\cNp$ for $\gR={\mbox{\tt
\fontsize{10pt}{0mm}\selectfont E}}_{8(8)}$}\\
\
\\
\begin{tabular}{c|c|c|c|c|c|c}

No. &${\rm Dyn}\,K\cdot x$
%%{\cal O}$
&${\rm Dyn}\,G\cdot x$

& $\dim^{}_{\C} K\cdot x$
%%{\cal O}$
& $\sharp (G\cdot x\cap \p)/K$ &
$\z_{\ka}(x)/{\rm rad}_u\z_{\ka}(x)$&
$\dim^{}_{\C}{\rm rad}_u\z_{\ka}(x)$
\\
\hline \hline
&&&&&\\[-13pt]
1&40000000&20000000&78&3&$2{\mbox{\tt
\fontsize{10pt}{0mm}\selectfont G}}_2$&14
\\
\hline
2&00000004&02000000&92&3&${\mbox{\tt
\fontsize{10pt}{0mm}\selectfont A}}_2$&20
\\
\hline
3&21010100&00010001&96&3&${\mbox{\tt
\fontsize{10pt}{0mm}\selectfont
A}}_1+{\mbox{\tt
\fontsize{10pt}{0mm}\selectfont T}}_1$&20
\\
\hline
4&00400000&00000200&97&2&$2{\mbox{\tt
\fontsize{10pt}{0mm}\selectfont A}}_1$&17
\\
\hline
5&40000040&02000002&99&3&${\mbox{\tt
\fontsize{10pt}{0mm}\selectfont A}}_2$&13
\\
\hline
6&20200200&00002000&104&3&0&16\\[-3pt]
7&00004000&00002000&104&3&0&16\\[-3pt]
8&02002002&00002000&104&3&0&16\\
\hline
9&40040000&20000200&105&2&$2{\mbox{\tt
\fontsize{10pt}{0mm}\selectfont A}}_1$&9
\\
\hline
10&02002022&00002002&107&3&${\mbox{\tt
\fontsize{10pt}{0mm}\selectfont T}}_1$&12\\[-3pt]
11&00400040&00002002&107&3&${\mbox{\tt
\fontsize{10pt}{0mm}\selectfont T}}_1$&12\\
\hline
12&31010211&10010101&108&2&${\mbox{\tt
\fontsize{10pt}{0mm}\selectfont T}}_1$&11\\
\hline
13&13111101&10010102&109&2&${\tt T}_1$&10\\
\hline
14&20202022&00020002&110&2&0&10\\[-3pt]
15&04004000&00020002&110&2&0&10\\
\hline
16&02022022&20002002&111&3&${\mbox{\tt
\fontsize{10pt}{0mm}\selectfont T}}_1$&8\\[-3pt]
17&40040040&20002002&111&3&${\mbox{\tt
\fontsize{10pt}{0mm}\selectfont T}}_1$&8
\end{tabular}
\end{center}
}

{\fontsize{9pt}{6mm}\selectfont
\begin{center}{TABLE 14 ({\it\!continued})} \\
{$(-1)$-distinguished $K$-orbits $
%%{\cal O}=
K\cdot x$ in $\cNp$ for $\gR={\mbox{\tt
\fontsize{10pt}{0mm}\selectfont E}}_{8(8)}$}\\
\
\\
\begin{tabular}{c|c|c|c|c|c|c}

No. &${\rm Dyn}\,K\cdot x$
%%{\cal O}$
&${\rm Dyn}\,G\cdot x$

& $\dim^{}_{\C} K\cdot x$
%%{\cal O}$
& $\sharp (G\cdot x\cap \p)/K$ &
$\z_{\ka}(x)/{\rm rad}_u\z_{\ka}(x)$&
$\dim^{}_{\C}{\rm rad}_u\z_{\ka}(x)$
\\
\hline \hline
&&&&&\\[-13pt]

18&00400400&00020020&112&2&0&8\\[-3pt]
19&22202022&00020020&112&2&0&8\\
\hline
20&22202042&00020022&113&2&0&7\\[-3pt]
21&04004040&00020022&113&2&0&7\\
\hline
22&22222022&20020020&114&2&0&6\\[-3pt]
23&40040400&20020020&114&2&0&6\\
\hline
24&22222042&20020022&115&2&0&5\\[-3pt]
25&04040044&20020022&115&2&0&5\\
\hline
26&22222222&20020202&116&2&0&4\\[-3pt]
27&44040400&20020202&116&2&0&4\\
\hline
28&24222242&20020222&117&2&0&3\\[-3pt]
29&44040440&20020222&117&2&0&3\\
\hline
30&44044044&22202022&118&1&0&2\\
\hline
31&44440444&22202222&119&1&0&1\\
\hline
32&84444444&22222222&120&1&0&0\\
\end{tabular}
\end{center}
}

\vskip 5mm

{\fontsize{9pt}{6mm}\selectfont
\begin{center}{TABLE 15 } \\
{$(-1)$-distinguished $K$-orbits $
%%{\cal O}=
K\cdot x$ in $\cNp$ for $\gR={\mbox{\tt
\fontsize{10pt}{0mm}\selectfont E}}_{8(-24)}$}\\
\
\\
\begin{tabular}{
c|c|c|c|c|c|c}

No. &${\rm Dyn}\,K\cdot x$
%%{\cal O}$
&${\rm Dyn}\, G\cdot x$ &
$\dim^{}_{\C}K\cdot x$
%%{\cal O}$
& $\sharp (G\cdot x\cap \p)/K$ &
$\z_{\ka}(x)/{\rm rad}_u\z_{\ka}(x)$&
$\dim^{}_{\C}{\rm rad}_u\z_{\ka}(x)$
\\
\hline \hline
&&&&&\\[-13pt]
1&00000004&00000002&57&3&${\mbox{\tt
\fontsize{10pt}{0mm}\selectfont E}}_6$&1
\\
\hline
2&00000204&00000020&83&3&${\mbox{\tt
\fontsize{10pt}{0mm}\selectfont D}}_4$&25
\\[-3pt]
3&00000040&00000020&83&3&${\mbox{\tt
\fontsize{10pt}{0mm}\selectfont D}}_4$&25
\\
\hline
4&00000048&00000022&84&2&${\mbox{\tt
\fontsize{10pt}{0mm}\selectfont F}}_4$&0
\\
\hline
5&01100012&10000100&89&1&${\mbox{\tt
\fontsize{10pt}{0mm}\selectfont
B}}_2+{\mbox{\tt
\fontsize{10pt}{0mm}\selectfont T}}_1$&36
\\
\hline
6&40000004&20000002&90&2&${\mbox{\tt
\fontsize{10pt}{0mm}\selectfont A}}_4$&22
\\
\hline
7&10100111&10000101&94&1&${\mbox{\tt
\fontsize{10pt}{0mm}\selectfont
A}}_2+{\mbox{\tt
\fontsize{10pt}{0mm}\selectfont T}}_1$&33
\\
\hline
8&01100034&10000102&95&1&${\mbox{\tt
\fontsize{10pt}{0mm}\selectfont A}}_3$&26
\\
\hline
9&00020000&00000200&97&1&$2{\mbox{\tt
\fontsize{10pt}{0mm}\selectfont
A}}_1$&33\\
\hline
10&20000222&20000020&99&2&${\mbox{\tt
\fontsize{10pt}{0mm}\selectfont G}}_2$&23
\\[-3pt]
11&00000404&20000020&99&2&${\mbox{\tt
\fontsize{10pt}{0mm}\selectfont G}}_2$&23
\\
\hline
12&20000244&20000022&100&2&${\mbox{\tt
\fontsize{10pt}{0mm}\selectfont B}}_3$&15
\\[-3pt]
13&40000048&20000022&100&2&${\mbox{\tt
\fontsize{10pt}{0mm}\selectfont B}}_3$&15

\end{tabular}
\end{center}
}

\vskip 5mm

{\fontsize{9pt}{6mm}\selectfont
\begin{center}{TABLE 15 ({\it\!continued})} \\
{$(-1)$-distinguished $K$-orbits $
%%{\cal O}=
K\cdot x$ in $\cNp$ for $\gR={\mbox{\tt
\fontsize{10pt}{0mm}\selectfont E}}_{8(-24)}$}\\
\
\\
\begin{tabular}{
c|c|c|c|c|c|c}

No. &${\rm Dyn}\,K\cdot x$
%%{\cal O}$
&${\rm Dyn}\, G\cdot x$ &
$\dim^{}_{\C}K\cdot x$
%%{\cal O}$
& $\sharp (G\cdot x\cap \p)/K$ &
$\z_{\ka}(x)/{\rm rad}_u\z_{\ka}(x)$&
$\dim^{}_{\C}{\rm rad}_u\z_{\ka}(x)$
\\
\hline \hline
&&&&&\\[-13pt]
14&00020200&20000200&105&1&$2{\mbox{\tt
\fontsize{10pt}{0mm}\selectfont A}}_1$&25
\\
\hline
15&40000404&20000202&107&1&${\mbox{\tt
\fontsize{10pt}{0mm}\selectfont A}}_2$&21
\\
\hline
16&40000448&20000222&108&1&${\mbox{\tt
\fontsize{10pt}{0mm}\selectfont G}}_2$&14
\\
\end{tabular}
\end{center}
}

\vskip 5mm

{\fontsize{9pt}{6mm}\selectfont
\begin{center}{TABLE 16} \\
{$(-1)$-distinguished $K$-orbits $
%%{\cal O}=
K \cdot x$ in $\cNp$ for $\gR={\mbox{\tt
\fontsize{10pt}{0mm}\selectfont F}}_{4(4)}$}\\
\
\\
\begin{tabular}{c|c|c|c|c|c|c}

No.&${\rm Dyn}\,K\cdot x$
%%{\cal O}$
&${\rm Dyn}\, G\cdot x$ &
$\dim^{}_{\C}K\cdot x$
%%{\cal O}$
& $\sharp (G\cdot x\cap \p)/K$ &
$\z_{\ka}(x)/{\rm rad}_u\z_{\ka}(x)$&
$\dim^{}_{\C}{\rm rad}_u\z_{\ka}(x)$
\\

\hline \hline
&&&&&\\[-13pt]

1&0004&2000&15&3&$ {\mbox{\tt
\fontsize{10pt}{0mm}\selectfont A}}_2$&1 \\

\hline

2&0040&0200&20&3&0&4\\[-3pt]
3&0204&0200&20&3&0&4\\[-3pt]
4&2022&0200&20&3&0&4\\

\hline

5&2200&0048&21&2&$ {\mbox{\tt
\fontsize{10pt}{0mm}\selectfont A}}_1$&0\\

\hline

6&0404&0202&22&2&0&2\\[-3pt]
7&2222&0202&22&2&0&2\\

\hline

8&2244&2202&23&2&0&1\\[-3pt]
9&4048&2202&23&2&0&1\\

\hline

10&4448&2222&24&1&0&0\\

\end{tabular}
\end{center}
}

\vskip 5mm
%%\newpage

{\fontsize{9pt}{6mm}\selectfont
\begin{center}{TABLE 17} \\
{ $(-1)$-distinguished $K$-orbits $
%%{\cal O}=
K \cdot x$ in $\cNp$ for $\gR={\mbox{\tt
\fontsize{10pt}{0mm}\selectfont F}}_{4(-20)}$}\\
\
\\
\begin{tabular}{c|c|c|c|c|c|c}

No.&${\rm Dyn}\,K\cdot x$
%%{\cal O}$
&${\rm Dyn}\, G\cdot x$ &
$\dim^{}_{\C}K\cdot x$
%%{\cal O}$
& $\sharp (G\cdot x\cap \p)/K$ &
$\z_{\ka}(x)/{\rm rad}_u\z_{\ka}(x)$&
$\dim^{}_{\C}{\rm rad}_u\z_{\ka}(x)$
\\

\hline \hline
&&&&&\\[-13pt]

1&0001&0001&11&1& ${\mbox{\tt
\fontsize{10pt}{0mm}\selectfont A}}_3$&10  \\

\hline

2&4000&0002&15&1& ${\mbox{\tt
\fontsize{10pt}{0mm}\selectfont G}}_2$ &7\\
\end{tabular}
\end{center}
}

 \vskip 5mm

{\fontsize{9pt}{6mm}\selectfont
\begin{center}{TABLE 18} \\
{$(-1)$-distinguished $K$-orbits $
%%{\cal
%%O}=
K\cdot x$ in $\cNp$ for $\gR={\mbox{\tt
\fontsize{10pt}{0mm}\selectfont G}}_{2(2)}$}\\
\
\\
\begin{tabular}{c|c|c|c|c|c|c}
No.&${\rm Dyn}\,K\cdot x$
%%{\cal O}$
&${\rm Dyn}\, G\cdot x$ &
$\dim^{}_{\C}K\cdot x$
%%{\cal O}$
& $\sharp (G\cdot x\cap \p)/K$ &
$\z_{\ka}(x)/{\rm rad}_u\z_{\ka}(x)$&
$\dim^{}_{\C}{\rm rad}_u\z_{\ka}(x)$\\
\hline \hline
&&&&& \\[-13pt]
1&22&02&5&2&0&1\\[-3pt]
2&04&02&5&2&0&1\\
\hline
3&48&22&6&1&0&0\\
\end{tabular}
\end{center}
}

{\fontsize{10pt}{2mm}\selectfont

}

\end{document}